\documentclass[12pt]{article}
\usepackage{full page, latexsym,
amscd,
amsthm,
graphicx, epsfig, amssymb}

\newtheorem{thm}{Theorem}[section]
\newtheorem{defn}[thm]{Definition}
\newtheorem{prop}[thm]{Proposition}
\newtheorem{cor}[thm]{Corollary}
\newtheorem{lemma}[thm]{Lemma}
\newtheorem{rema}[thm]{Remark}

\newtheorem{exam}[thm]{Example}

\newcommand{\halmos}{\rule{1ex}{1.4ex}}

\newcommand{\nn}{\nonumber \\}

 \newcommand{\res}{\mbox{\rm Res}}

\renewcommand{\hom}{\mbox{\rm Hom}}

 \newcommand{\pf}{{\it Proof.}\hspace{2ex}}
 \newcommand{\epf}{\hspace*{\fill}\mbox{$\halmos$}}
 \newcommand{\epfv}{\hspace*{\fill}\mbox{$\halmos$}\vspace{1em}}

\newcommand{\wt}{\mbox{\rm wt}\ }
\newcommand{\swt}{\mbox{\rm {\scriptsize wt}}\ }

\newcommand{\lbar}{\bigg\vert}

\newcommand{\C}{\mathbb{C}}
\newcommand{\Z}{\mathbb{Z}}
\newcommand{\R}{\mathbb{R}}

\newcommand{\one}{\mathbf{1}}

\title{ {\bf Two constructions of grading-restricted vertex (super)algebras} }
\date{}
\author{Yi-Zhi Huang}

\begin{document}

\bibliographystyle{alpha}
\maketitle
\begin{abstract}
We give two constructions of grading-restricted vertex (super)algebras. We first give
a new construction of a class of grading-restricted vertex (super)algebras
originally obtained by Meurman and Primc using a different method.
This construction is based on a new definition of vertex operators
and a  new method.
Our second construction is a generalization
of the author's construction of the moonshine module vertex operator algebra
and a related vertex operator superalgebra. This construction needs
properties of intertwining operators formulated and proved by the author.
\end{abstract}

\renewcommand{\theequation}{\thesection.\arabic{equation}}
\renewcommand{\thethm}{\thesection.\arabic{thm}}
\setcounter{equation}{0}
\setcounter{thm}{0}
\section{Introduction}

Vertex operator (super)algebras are algebraic structures formed by
meromorphic fields in two-dimensional (super)conformal field theories.
Based on this connection, a program to construct
two-dimensional (super)conformal field theories
using the representation theory of vertex operator (super)algebras
has been successfully carried out
during the past 25 years, though there are still many open problems to be solved.
For a nontechnical  brief description of this program, see
the author's blog article \cite{Hblog}. The construction
of vertex operator algebras is a prerequisite for this program.

A vertex operator (super)algebra is a grading-restricted vertex (super)algebra equipped with
a conformal element. In this paper, we are interested only in the construction
of vertex operators, not the construction of conformal elements. Thus we shall restrict
our attention to grading-restricted vertex (super)algebras instead of vertex operator (super)algebras.

The existing constructions of grading-restricted vertex (super)algebras can mostly be divided
into three types. The first type includes the constructions of grading-restricted vertex
(super)algebras associated to the Heisenberg algebras, Clifford algebras
and lattices (see \cite{FLM} and \cite{FFR}). These constructions are based on
explicit formulas for vertex operators. The second type includes the constructions
of the grading-restricted vertex (super)algebras associated to
the Virasoro algebra, affine Lie algebras and superconformal algebras (see \cite{FZ}, \cite{KW} and \cite{A}).
These constructions
are based on a definition of vertex operators using generating fields and the residue of $x_{1}$
of the Jacobi identity for vertex (super)algebras. The third type is the so-called orbifold
construction,  including
the constructions of the moonshine module vertex operator algebra
and a related vertex operator superalgebra (see \cite{FLM}, \cite{DGM}
\cite{H1} and \cite{H2}). These constructions
are based on suitable intertwining operators among fixed-point vertex operator
subalgebras and modules for the subalgebras.

In this paper, we first give a new construction of a class of grading-restricted vertex (super)algebras. Our construction
is based on a new definition of vertex operators using the rational functions
obtained from the products of generating fields. This definition is motivated by the
associativity for vertex operators. Our construction is also based on a new method.
A uniqueness result
shows that the class of grading-restricted vertex algebras given by this construction is the same
as the class given by Meurman and Primc \cite{MP}  using a different construction of the second type.

We also give  in this paper a generalization of the author's construction in \cite{H1} of the
moonshine module vertex operator algebra. Given a grading-restricted
vertex algebra equipped with a compatible $\mathfrak{sl}(2, \C)$ action
(called quasi-vertex operator algebra  in \cite{FHL})
and a nondegenerate symmetric invariant bilinear form
and a module equipped with similar additional structures
satisfying certain conditions, we construct a quasi-vertex operator algebra or superalgebra
structure on the direct sum of  the algebra and the module using
suitable intertwining operators. Examples of such quasi-vertex operator algebras or superalgebras
include, as mentioned above, the moonshine module vertex operator algebra and a related
vertex operator superalgebra. The vertex operator
superalgebras associated to the Moore-Read states \cite{MR} in the study of quantum Hall states
can also be obtained using this construction.

In addition to being useful for obtaining examples of grading-restricted vertex (super)algebras,
quasi-vertex operator (super)algebras and
vertex operator (super)algebras, these two constructions provide new methods in the study of
these algebras and their
representations. For example, the author's cohomology theory for grading-restricted
vertex algebras \cite{H8} and \cite{H9} are constructed and developed based on the rational functions of
the products and iterates of vertex operators and operators corresponding to
cochains. Our first construction will  be useful in the study of this cohomology theory for
those algebras that can be constructed in this way.
The (logarithmic) tensor category theory for module categories for vertex operator algebras (or more generally
for quasi-vertex operator algebras or M\"{o}bius vertex algebras) is based on (logarithmic) intertwining operators
(see the expositions \cite{HL} and \cite{HLZ} and the references in these papers for details on this theory).
Our second construction uses intertwining operators and is indeed closely related to
the tensor category theory. The second construction can be generalized to more complicated
orbifold constructions and  these generalizations will be studied in future publications.

In the present paper, instead of the formal variable approach, we use the complex analysis
approach. Although we do use residues, our method depends mainly on
the method of analytic
extensions and the properties of rational functions of the special type appearing
in the theory of vertex operator algebras.  The advantage of our analytic  approach is that
every definition or proof has a geometric meaning. The geometric meaning is not logically
needed in this paper and in many papers studying vertex operator algebras. But it often
provides ideas and motivations for many constructions and proofs. One
important example is the
cohomology theory for grading-restricted vertex algebras \cite{H8} \cite{H9} mentioned above
and the author's ongoing study of the representation theory of vertex operator algebras
using this cohomology theory. The geometric insights play a crucial role in obtaining
new results and developing new methods in this theory.

The present paper is organized as follows: In the next section, we recall the
definitions of grading-restricted vertex (super)algebra, module and intertwining operator
and some other useful notions. In Sections 3 and 4, we give our first and second constructions,
respectively.

\renewcommand{\theequation}{\thesection.\arabic{equation}}
\renewcommand{\thethm}{\thesection.\arabic{thm}}
\setcounter{equation}{0}
\setcounter{thm}{0}
\section{Grading-restricted vertex (super)algebras}

In this section, we recall the definitions of grading-restricted
vertex (super)algebra, module and intertwining operator. We also
recall some other notions, including the notions of quasi-vertex operator (super)algebra,
conformal element, vertex operator (super)algebra,
fusion rule and nondegenerate symmetric invariant bilinear form.
Except for conformal elements and vertex operator algebras, these are all needed in our constructions.

For a $\mathbb{Z}$-graded vector space $V=\coprod_{n\in \Z}V_{(n)}$, as usual, we
denote its graded dual space $\coprod_{n\in \Z}V_{(n)}^{*}$ by $V'$ and its
algebraic completion $\prod_{n\in \Z}V_{(n)}$ by $\overline{V}$. On
$V$ and $V'$, we use the topology given by the dual pair $(V, V')$. In the case
that $\dim V_{(n)}<\infty$ for $n\in \Z$,
we use the topology on $\overline{V}=(V')^{*}$ given by the dual pair $(V', \overline{V})$.
We give a topology to the space $\hom(V^{\otimes n}, \overline{V})$ by identifying it with
a subspace of the dual space of $V^{\otimes n}\otimes V'$ and then using the topology
induced from the dual pair of this subspace and $V^{\otimes n}\otimes V'$. More specifically,
a sequence (or more generally a net)
$\{f_{n}\}$ in $\hom(V\otimes \cdots \otimes V, \overline{V})$ is convergent
to $f\in \hom(V\otimes \cdots \otimes V, \overline{V})$ if for $v_{1}, \dots,
v_{n}\in V$ and $v'\in V'$, $\langle v', f_{n}(v_{1}\otimes \cdots\otimes v_{n})\rangle$
is convergent to $\langle v', f(v_{1}\otimes \cdots\otimes v_{n})\rangle$.
In particular, analytic maps from a region in $\C$ to $\hom(V^{\otimes n}, \overline{V})$
make sense. For a $\C$-graded vector space, we use
the same notations and topologies.

We give the definition of grading-restricted vertex (super)algebra first.

\begin{defn}\label{grvsa}
{\rm A {\it grading-restricted vertex superalgebra} is
a $\frac{\Z}{2}$-graded vector space
$V=\coprod_{n\in \frac{\Z}{2}}V_{(n)}=V^{0}\oplus V^{1}$, where $V^{0}=\coprod_{n\in \Z}V_{(n)}$
and $V^{1}=\coprod_{n\in \Z+\frac{1}{2}}V_{(n)}$,
equipped with
an analytic map
\begin{eqnarray*}
Y_{V}: \C^{\times}&\to &\hom(V\otimes V, \overline{V}), \nn
z&\mapsto &Y_{V}(\cdot, z)\cdot: u\otimes v\mapsto Y_{V}(u, z)v
\end{eqnarray*}
called {\it vertex operator map}
and a {\it vacuum} $\mathbf{1}\in V_{(0)}$
satisfying the following axioms:

\begin{enumerate}

\item Axioms for the grading:
(a) {\it Grading-restriction condition}: When $n$ is sufficiently negative,
$V_{(n)}=0$ and $\dim V_{(n)}<\infty$ for $n\in \frac{\Z}{2}$.
(b) {\it $L(0)$-bracket formula}: Let $L_{V}(0): V\to V$
be defined by $L_{V}(0)v=nv$ for $v\in V_{(n)}$. Then
$$[L_{V}(0), Y_{V}(v, z)]=\frac{d}{dz}Y_{V}(v, z)+Y_{V}(L_{V}(0)v, z)$$
for $v\in V$.


\item Axioms for the vacuum: (a) {\it Identity property}:
Let $1_{V}$ be the identity operator on $V$. Then
$Y_{V}(\mathbf{1}, z)=1_{V}$. (b)
{\it Creation property}: For $u\in V$,
$\lim_{z\to 0}Y_{V}(u, z)\mathbf{1}$ exists and is equal to $u$.

\item {\it $L(-1)$-derivative property}:
Let $L_{V}(-1): V\to V$ be the operator
given by
$$L_{V}(-1)v=\lim_{z\to 0}\frac{d}{dz}Y_{V}(v, z)\one$$
for $v\in V$. Then for $v\in V$,
$$\frac{d}{dz}Y_{V}(v, z)=Y_{V}(L_{V}(-1)v, z)=[L_{V}(-1), Y_{V}(v, z)].$$

\item {\it Duality}: For $u_{1}\in V^{|u_{1}|}$, $u_{2}\in V^{|u_{2}|}$, and $v\in V$ where $|u_{1}|, |u_{2}|$  are $0$ or $1$,
$v'\in V'$, the series
\begin{eqnarray*}
&\langle v', Y_{V}(u_{1}, z_{1})Y_{V}(u_{2}, z_{2})v\rangle,&\\
&
(-1)^{|u_{1}| |u_{2}|}\langle v', Y_{V}(u_{2}, z_{2})Y_{V}(u_{1}, z_{1})v\rangle,&\\
&\langle v', Y_{V}(Y_{V}(u_{1}, z_{1}-z_{2})u_{2}, z_{2})v\rangle&
\end{eqnarray*}
are absolutely convergent
in the regions $|z_{1}|>|z_{2}|>0$, $|z_{2}|>|z_{1}|>0$,
$|z_{2}|>|z_{1}-z_{2}|>0$, respectively, to a common rational function
in $z_{1}$ and $z_{2}$ with the only possible poles at $z_{1}, z_{2}=0$ and
$z_{1}=z_{2}$.

\end{enumerate}
In the case that $V^{1}=0$, the grading-restricted vertex superalgebra just defined is called a {\it grading-restricted
vertex algebra}.}
\end{defn}

We denote the grading-restricted  vertex (super)algebra
just defined by $(V, Y_{V},  \one)$ or by $V$. Note that in the definition above,
we use the duality instead of the Jacobi identity or weak commutativity as the main
axiom.

Although we are mainly interested in grading-restricted vertex (super)algebras in this paper,
in our second construction, we need an operator $L_{V}(1)$ acting on the
algebra $V$ such that together with $L_{V}(0)$ and $L_{V}(-1)$, it gives an action of the Lie
algebra $\mathfrak{sl}(2, \C)$ on $V$ and satisfies the usual bracket formula between the basis of
$\mathfrak{sl}(2, \C)$ and vertex operators.
Also our main examples all have conformal elements. Because of these reasons, we also recall the definitions
of quasi-vertex operator (super)algebra, conformal element and vertex operator (super)algebra from
\cite{FHL}.

\begin{defn}\label{qvsa}
{\rm A {\it quasi-vertex operator (super)algebra} is a grading-restricted vertex (super)algebra
$(V, Y_{V},  \one)$ together with an operator $L_{V}(1)$ of weight $1$ on $V$ satisfying
\begin{eqnarray*}
[L_{V}(-1), L_{V}(1)]&=&-2L_{V}(0),\\
{[L_{V}(1), Y_{V} (v, z)]}&= &Y_{V} (L_{V}(1)v, z) + 2zY_{V} (L_{V}(0)v, z) + z^{2}Y_{V} (L_{V}(-1)v, z)
\end{eqnarray*}
for $v\in V$. }
\end{defn}

We denote the  quasi-vertex operator (super)algebra
just defined by $(V, Y_{V},  \one, L_{V}(1))$ or simply by $V$.

\begin{defn}\label{vosa}
{\rm Let $(V, Y_{V},  \one)$ be a grading-restricted vertex (super)algebra.
A {\it conformal element} of $V$ is an element $\omega\in V$
satisfying the following axioms:

\begin{enumerate}

\item There exists $c\in \C$ such that
$Y (\omega, z)\omega$ expanded as a $V$-valued Laurent series is equal
to $L_{V}(-1)\omega z^{-1} + 2\omega z^{-2} + \frac{c}{2}\one z^{-4}$
plus a $V$-valued power series in $z$.

\item  $L_{V}(-1) = \res_{z}Y_{V}(\omega, z)$ and $L_{V}(0)=
\res_{z}z Y_{V}(\omega, z)$ ($\res_{z}$ being the operation of taking the coefficient of $z^{-1}$
of a Laurent series).

\end{enumerate}
A grading-restricted
vertex (super)algebra equipped with a conformal element is called a {\it vertex operator
(super)algebra} (or, more consistently, {\it grading-restricted conformal vertex (super)algebra}). }
\end{defn}

We denote the  vertex operator (super)algebra
just defined by $(V, Y_{V},  \one, \omega)$ or simply by $V$.

\begin{rema}\label{skew-symm}
{\rm The absolute convergence of
$$\langle v', Y_{V}(u_{1}, z_{1})Y_{V}(u_{2}, z_{2})v\rangle$$
and
$$\langle v', Y_{V}(Y_{V}(u_{1}, z_{1}-z_{2})u_{2}, z_{2})v\rangle$$
to rational functions in the regions $|z_{1}|>|z_{2}|>0$ and $|z_{2}|>|z_{1}-z_{2}|>0$, respectively,
in Definition \ref{grvsa}  are called the {\it rationality of
the products} and  the {\it rationality of
the iterates}, respectively. The statement that
$$\langle v', Y_{V}(u_{1}, z_{1})Y_{V}(u_{2}, z_{2})v\rangle$$
and
$$(-1)^{|u_{1}| |u_{2}|}\langle v', Y_{V}(u_{2}, z_{2})Y_{V}(u_{1}, z_{1})v\rangle$$
converges  in the regions $|z_{1}|>|z_{2}|>0$ and $|z_{2}|>|z_{1}|>0$, respectively, to a common rational function
is called the {\it commutativity}. The statement that
$$\langle v', Y_{V}(u_{1}, z_{1})Y_{V}(u_{2}, z_{2})v\rangle$$
and
$$\langle v', Y_{V}(Y_{V}(u_{1}, z_{1}-z_{2})u_{2}, z_{2})v\rangle$$
converges  in the regions $|z_{1}|>|z_{2}|>0$ and $|z_{2}|>|z_{1}-z_{2}|>0$, respectively, to a common rational function
is called the {\it associativity}. In fact, with all the other properties still hold, it is easy to
show that
the commutativity is equivalent to the {\it skew-symmetry}: For $u\in V^{|u|}$ and $v\in V^{|v|}$,
$$Y_{V}(u, z)v=(-1)^{|u||v|}e^{zL_{V}(-1)}Y_{V}(v, -z)u.$$
In particular, we can replace the duality in Definition
\ref{grvsa} by the rationality of the products and iterates, the associativity and the skew-symmetry.
We shall need this fact below.}
\end{rema}

Next we give the definition of module for a grading-restricted vertex (super)algebra.

\begin{defn}\label{module}
{\rm Let $V$ be a grading-restricted vertex superalgebra.
A {\it $V$-module} is
a $\C\times \Z_{2}$-graded vector space
$$W=\coprod_{n\in \C, \alpha\in \Z_{2}}W^{\alpha}_{(n)}=
\coprod_{n\in \C}W_{(n)}=W^{0}\oplus W^{1}$$
(where $W_{(n)}=W_{(n)}^{0}\oplus W_{(n)}^{1}$, $W^{0}=\coprod_{n\in \C}W^{0}_{(n)}$
and $W^{1}=\coprod_{n\in \C}V^{1}_{(n)}$)
equipped with
a {\it vertex operator map}
\begin{eqnarray*}
Y_{W}: \C^{\times} &\to& \hom(V\otimes W, \overline{W}), \nn
z&\mapsto & Y_{W}(\cdot, z)\cdot: u\otimes w\mapsto Y_{W}(u, z)w
\end{eqnarray*}
satisfying the following axioms:

\begin{enumerate}

\item Axioms for the gradings:
(a) {\it Grading-restriction condition}: When the real part of $n$ is sufficiently negative,
$W_{(n)}=0$ and $\dim W_{(n)}<\infty$ for $n\in \C$. (b) {\it $L(0)$-bracket formula}: Let $L_{W}(0): V\to V$
be defined by $L_{W}(0)v=nv$ for $v\in W_{(n)}$. Then
$$[L_{W}(0), Y_{W}(v, z)]=\frac{d}{dz}Y_{W}(v, z)+Y_{W}(L_{V}(0)v, z)$$
for $v\in V$. (c) {\it Grading compatibility}: For $\alpha, \beta\in \Z_{2}$, $u\in V^{\alpha}$
and $w\in W^{\beta}$, $Y_{W}(u, z)w\in \overline{W^{\alpha+\beta}}$.


\item  {\it Identity property}:
Let $1_{W}$ be the identity operator on $W$. Then
$Y_{W}(\mathbf{1}, z)=1_{W}$.

\item {\it $L(-1)$-derivative property}:
There exists $L_{W}(-1): W\to W$ such that for $u\in V$,
$$\frac{d}{dz}Y_{W}(u, z)=Y_{W}(L_{V}(-1)u, z)=[L_{W}(-1), Y_{W}(u, z)].$$

\item {\it Duality}: For $u_{1}\in V^{|u_{1}|}$, $u_{2}\in V^{|u_{2}|}$, and $w\in W$ where $|u_{1}|, |u_{2}|$  are $0$ or $1$,
$w'\in W'$, the series
\begin{eqnarray*}
&\langle w', Y_{W}(u_{1}, z_{1})Y_{W}(u_{2}, z_{2})w\rangle,&\\
&
(-1)^{|u_{1}| |u_{2}|}\langle w', Y_{W}(u_{2}, z_{2})Y_{W}(u_{1}, z_{1})w\rangle,&\\
&\langle w', Y_{W}(Y_{V}(u_{1}, z_{1}-z_{2})u_{2}, z_{2})w\rangle&
\end{eqnarray*}
are absolutely convergent
in the regions $|z_{1}|>|z_{2}|>0$, $|z_{2}|>|z_{1}|>0$,
$|z_{2}|>|z_{1}-z_{2}|>0$, respectively, to a common rational function
in $z_{1}$ and $z_{2}$ with the only possible poles at $z_{1}, z_{2}=0$ and
$z_{1}=z_{2}$.

\end{enumerate}
When $V$ is a grading-restricted vertex algebra, a {\it $V$-module} is a $V$-module
$W$  with $W^{1}=0$ when $V$ is viewed as a grading-restricted vertex superalgebra.
When $V$ is a vertex operator (super)algebra, a {\it $V$-module}
is a $V$-module when $V$ is viewed as a grading-restricted vertex (super)algebra.
When $V$ is a quasi-vertex operator (super)algebra, a {\it $V$-module} is a $V$-module $W$
when $V$ is viewed as a grading-restricted vertex (super)algebra together with an
operator $L_{W}(1)$ of weight $1$ on $W$ satisfying
\begin{eqnarray*}
[L_{W}(-1), L_{W}(1)]&=&-2L_{W}(0),\\
{[L_{W}(1), Y_{W} (v, z)]}&= &Y_{W} (L_{V}(1)v, z) + 2zY_{W} (L_{V}(0)v, z) + z^{2}Y_{W} (L_{V}(-1)v, z)
\end{eqnarray*}
for $v\in V$.}
\end{defn}

We denote the $V$-module just defined by $(W, Y_{W})$ or simply by $W$.
In the case that $V$ is a quasi-vertex operator (super)algebra, we denote the
$V$-module just defined by $(W, Y_{W}, L_{W}(1))$ or simply by $W$.

\begin{rema}
{\rm In Definition \ref{module}, as in Definition \ref{grvsa}, we can also separate the rationality of products and iterates, the
commutativity and the associativity (see Remark \ref{skew-symm}). In fact, for modules, it is easy to see that the commutativity is a consequence
of the rationality, the associativity and the skew-symmetry for $V$. Since the skew-symmetry for $V$
always holds, the duality in Definition \ref{module}
can be replaced by the rationality and the associativity. }
\end{rema}

We also need the important notion of intertwining operator. We could define an intertwining
operator to be an analytic map from $\C^{\times}$ to a suitable space of linear maps. But
in general it is not single valued and we need to choose a preferred branch. Because of
this complication, we shall define an intertwining operator as usual to be a linear
map to a space of formal series with complex powers. After the definition, we shall
explain how to choose a special branch of the intertwining operator to obtain an analytic
map from $\C^{\times}$ to the corresponding space of linear maps. These maps are what
we shall use in this paper.

\begin{defn}\label{int-op}
{\rm Let $V$ be a grading-restricted vertex superalgebra (a grading-restricted vertex algebra being
a special case) and $W_{1},
W_{2}, W_{3}$ $V$-modules.
An {\it intertwining operator of type ${W_{3}\choose W_{1}W_{2}}$} is
a linear map
\begin{eqnarray*}
\mathcal{Y}: W_{1}\otimes W_{2} &\to  & W_{3}\{x\} \nn
w_{1}\otimes w_{2} & \mapsto  &\mathcal{Y}(w_{1}, x)w_{2}
\end{eqnarray*}
(where $W_{3}\{x\}$ is the space of formal series of the form
$\sum_{n\in \C}a_{n}x^{n}$ for $a_{n}\in W_{3}$ and $x$
is a formal variable)
satisfying the following axioms:

\begin{enumerate}

\item  {\it $L(0)$-bracket formula}: For $w_{1}\in W_{1}$,
$$L_{W_{3}}(0)\mathcal{Y}(w_{1}, x)-\mathcal{Y}(w_{1}, x)L_{W_{2}}(0)
=\frac{d}{dx}\mathcal{Y}(w_{1}, x)+Y_{W}(L_{W_{1}}(0)w_{1}, x).$$


\item {\it $L(-1)$-derivative property}:
For $w_{1}\in W_{1}$,
$$\frac{d}{dx}\mathcal{Y}(w_{1}, x)=Y_{W}(L_{W_{1}}(-1)w_{1}, x)=L_{W_{3}}(-1) \mathcal{Y}(w_{1}, x)-
\mathcal{Y}(w_{1}, x)L_{W_{2}}(-1) .$$

\item {\it Duality with vertex operators}: For $u\in V^{|u|}$, $w_{1}\in W^{|w_{1}|}$, and $w_{2}\in W$
where $|u|, |w_{1}|$  are $0$ or $1$,
$w_{3}'\in W_{3}'$, for any single-valued branch $l(z_{2})$ of the logarithm of $z_{2}$ in the region
$z_{2}\ne 0$, $0\le \arg z_{2}\le 2\pi$,
the series
\begin{eqnarray*}
&\langle w_{3}', Y_{W_{3}}(u, z_{1})\mathcal{Y}(w_{1}, x_{2})w_{2}\rangle
\lbar_{x_{2}^{n}=e^{nl(z_{2})}, n\in \C},&\\
&(-1)^{|u| |w_{1}|}\langle w_{3}', \mathcal{Y}(w_{1}, x_{2})Y_{W_{2}}(u, z_{1})w_{2}\rangle
\lbar_{x_{2}^{n}=e^{nl(z_{2})}, n\in \C},&\\
&\langle w_{3}', \mathcal{Y}(Y_{W_{1}}(u, z_{1}-z_{2})w_{1}, x_{2})w_{2}\rangle
\lbar_{x_{2}^{n}=e^{nl(z_{2})}, n\in \C}&
\end{eqnarray*}
are absolutely convergent
in the regions $|z_{1}|>|z_{2}|>0$, $|z_{2}|>|z_{1}|>0$,
$|z_{2}|>|z_{1}-z_{2}|>0$, respectively, to a common analytic function in $z_{1}$ and $z_{2}$
and can be analytically extended to a
multivalued analytic functions with the only possible poles $z_{1}=0$ and
$z_{1}=z_{2}$ and the only possible branch point
$z_{2}=0$.

\end{enumerate}
When $V$ is a quasi-vertex operator superalgebra (a quasi-vertex operator superalgebra being
a special case) and $W_{1}, W_{2}, W_{3}$ are $V$-modules,
an {\it intertwining operator of type ${W_{3}\choose W_{1}W_{2}}$} is an intertwining
operator $\mathcal{Y}$ of the same type in the sense above satisfying in addition the following
condition:
\begin{enumerate}
\setcounter{enumi}{3}

\item {\it $L(1)$-bracket formula}: For $w_{1}\in W_{1}$,
\begin{eqnarray*}
\lefteqn{L_{W_{3}}(1)\mathcal{Y} (w_{1}, z)-\mathcal{Y} (w_{1}, z)L_{W_{2}}(1)}\nn
&&= \mathcal{Y} (L_{W_{1}}(1)w_{1}, z) + 2z\mathcal{Y} (L_{W_{1}}(0)w_{1}, z) + z^{2}\mathcal{Y} (L_{W_{1}}(-1)w_{1}, z).
\end{eqnarray*}

\end{enumerate}
The dimension of the space of all intertwining operators of type ${W_{3}\choose W_{1}W_{2}}$
is called the {\it fusion rule of type ${W_{3}\choose W_{1}W_{2}}$} and is denoted
by $N_{W_{1}W_{2}}^{W_{3}}$.}
\end{defn}

We shall need later a formula equivalent to the $L(0)$-bracket formula called
{\it $L(0)$-conjugation formula} for an intertwining operator $\mathcal{Y}$: For $a\in \C$,
$$e^{aL_{W_{3}}(0)}\mathcal{Y}(w_{1}, x)e^{-aL_{W_{2}}(0)}=\mathcal{Y}(e^{aL_{W_{1}}(0)}w_{1}, e^{a}x).$$

For $z\in \C^{\times}$, let $\log z=\log |z|+i \arg z$, where $0\le \arg z< 2\pi$.
Let $\mathcal{Y}$ be an intertwining operator of type ${W_{3}\choose W_{1}W_{2}}$.
Then for $z\in \C^{\times}$, $w_{1}\in W_{1}$ and $w_{2}\in W_{2}$,
we use $\mathcal{Y}(w_{1}, z)w_{2}$ to denote $\mathcal{Y}(w_{1}, x)w_{2}
\lbar_{x^{n}=e^{n\log z}, n\in \C}$. In particular, we have a map
from $\C^{\times}$ to $\hom(W_{1}\otimes W_{2}, \overline{W}_{3})$ given by
$z\mapsto \mathcal{Y}(w_{1}, z)w_{2}$. Using this notation,
the three expressions in the duality axiom for intertwining
operators can be written as
\begin{eqnarray*}
&\langle w_{3}', Y_{W_{3}}(u, z_{1})\mathcal{Y}(w_{1}, z_{2})w_{2}\rangle,&\\
&(-1)^{|u| |w_{1}|}\langle w_{3}', \mathcal{Y}(w_{1}, z_{2})Y_{W_{2}}(u, z_{1})w_{2}\rangle,&\\
&\langle w_{3}', \mathcal{Y}(Y_{W_{1}}(u, z_{1}-z_{2})w_{1}, z_{2})w_{2}\rangle.&
\end{eqnarray*}
We shall always use this notation in this paper.

Finally we define nondegenerate symmetric invariant bilinear forms
on modules for a quasi-vertex operator algebra.

\begin{defn}
{\rm Let $V$ be a quasi-vertex operator algebra and $W$ a $V$-module. A
{\it nondegenerate symmetric invariant bilinear form
on $W$} is a nondegenerate symmetric bilinear form $(\cdot, \cdot)_{W}: W\otimes W\to \C$
 on $W$ such that for $m, n\in \C$, $m\ne n$,
$w_{1}\in W_{(m)}$ and $w_{2}\in W_{(n)}$,   $(w_{1}, w_{2})_{W}=0$ and
for $v\in V$, $w_{1}, w_{2}\in W$,
\begin{eqnarray*}
(w_{1}, L_{W}(1)w_{2})_{W}&=&(L_{W}(-1)w_{1}, w_{2})_{W},\\
(w_{1}, Y_{W}(v, z)w_{2})_{W}&=&(Y_{W}(e^{zL_{V}(1)}(-z^{-2})^{L_{V}(0)}v, z^{-1})w_{1}, w_{2})_{W}.
\end{eqnarray*}
In particular, a {\it  nondegenerate symmetric invariant bilinear form
on $V$} is defined to be  a nondegenerate symmetric invariant bilinear form
on $V$ when $V$ is viewed as a $V$-module.}
\end{defn}

\renewcommand{\theequation}{\thesection.\arabic{equation}}
\renewcommand{\thethm}{\thesection.\arabic{thm}}
\setcounter{equation}{0}
\setcounter{thm}{0}
\section{The first construction}

In this section, we give our first construction of grading-restricted vertex (super)algebras.
By a uniqueness result, our construction gives the same class of grading-restricted vertex algebras
as in \cite{MP}. What is new in this section is the construction,  including
a new formula for the vertex operator map and a  new method.
As mentioned in the introduction, we use the complex analysis approach.
But since in this section,
we work only with rational functions, the results and the proofs in this section
still work for grading-restricted vertex (super)algebras over any field of characteristic $0$.

Let $V=\coprod_{n\in \frac{\Z}{2}}$ be a $\frac{\Z}{2}$-graded vector space such that $V_{(n)}=0$ for $n$
sufficiently negative and $\dim V_{(n)}<\infty$ for $n\in \frac{\Z}{2}$.
Since  $\dim V_{(n)}<\infty$ for $n\in \frac{\Z}{2}$, we have
$\overline{V}=(V')^{*}$. Elements of $V_{(n)}$
is said to have {\it weight $n$}. Elements of $V^{0}=\coprod_{n\in \Z}V_{(n)}$ are said to be {\it even} and
elements of $V^{1}=\coprod_{n\in \Z+\frac{1}{2}}V_{(n)}$ are said to be {\it odd}.
Let $L_{V}(0): V\to V$ be the operator defined by the grading on $V$, that is,
by $L_{V}(0)v=nv$ for $v\in V_{(n)}$. Then for $a\in \C$, the operator $e^{a L_{V}(0)}$ on $V$
defined by $e^{a L_{V}(0)}v=e^{a n}v$ for $v\in V_{(n)}$
has a natural extension to $\overline{V}$. For $n\in \Z$, we use $\pi_{n}$ to denote the
projection from $V$ or $\overline{V}$ to $V_{(n)}$.


An operator $O$ on $V$ satisfying $[L_{V}(0), O]=nO$ is said to have {\it weight $n$}. Similarly for operators on
the graded dual $V'$ of $V$.

\begin{lemma}\label{laurent}
Let $\phi$ be an analytic  map from $\C^{\times}$
to $\hom(V, \overline{V})$. If there exists $\wt \phi\in \Z$ such that
$$[L_{V}(0), \phi(z)]=z\frac{d}{dz}\phi(z)+(\wt \phi)\phi(z),$$
then we have a Laurent expansion
$\phi(z)=\sum_{n\in \Z}\phi_{n}z^{-n-1}$ where for $n\in \Z$,
$\phi_{n}\in \hom(V, V)$ is homogeneous of weight $\wt \phi-n-1$. Moreover, for $v\in V$,
$\phi(z)v$ as a Laurent series in $z$ has only finitely many negative power terms and for $v'\in V'$,
$\langle v', \phi(z)\cdot\rangle$ as a Laurent series with coefficients in $V'$ has only finitely many
 positive powers of $z$.
\end{lemma}
\pf
From the bracket formula for $L_{V}(0)$ and $\phi(z)$, we obtain that for
$a\in \C$,
$$e^{aL_{V}(0)}\phi(z)e^{-aL_{V}(0)}=e^{a(\swt \phi)}\phi(e^{a}z).$$
In particular, taking $a$ such that
$e^{-a}=z$, we have
$\phi(z)=e^{a(\swt \phi)}e^{-aL_{V}(0)}\phi(1)e^{aL_{V}(0)}$.  Let $\phi_{n}: V\to V$
be defined by $\phi_{n}v=\pi_{(\swt \phi)-n-1+m}\phi(1)v$ for $v\in V_{(m)}$.
Then $\phi_{n}$ is of weight $\wt \phi-n-1$ and $\sum_{n\in \Z}\phi_{n}=\phi(1)$.
Moreover,
$e^{a(\swt \phi)}e^{-aL_{V}(0)}\phi_{n}e^{aL_{V}(0)}v=\phi_{n}v z^{-n-1}$ for $n\in \Z$ and $v\in V$.
Thus
\begin{eqnarray*}
\phi(z)v&=&e^{a(\swt \phi)}e^{-aL_{V}(0)}\phi(1)e^{aL_{V}(0)}v\nn
&=&\sum_{n\in \Z}e^{a(\swt \phi)}e^{-aL_{V}(0)}\phi_{n}e^{aL_{V}(0)}v\nn
&=&\sum_{n\in \Z}\phi_{n}v z^{-n-1}
\end{eqnarray*}
for $v\in V$.

Since $V_{(n)}=0$ for $n$ sufficiently negative and the weight of $\phi_{n}$ is $\wt \phi-n-1$,
for $v\in V$, $\phi(z)v$ has only finitely many negative power terms and for $v'\in V'$,
$\langle v', \phi(z)\cdot\rangle$ as a Laurent series with coefficients in $V'$ has only finitely many
 positive powers of $z$.
\epfv

Let $\phi^{i}$ for $i\in I$ be analytic maps from $\C^{\times}$
to $\hom(V, \overline{V})$ and let $\one\in V_{(0)}$.
Assume that $V$, $\phi^{i}$ for $i\in I$ and $\one\in V_{(0)}$
satisfy the following conditions:

\begin{enumerate}

\item  For $i\in I$, there exists $\wt \phi^{i}\in \Z$ such that
$[L_{V}(0), \phi^{i}(z)]=z\frac{d}{dz}\phi^{i}(z)+(\wt \phi^{i})\phi^{i}(z)$.

\item There exists an operator $L_{V}(-1)$ on $V$ such that $L_{V}(-1)\one=0$ and
$[L_{V}(-1), \phi^{i}(z)]=\frac{d}{dz}\phi^{i}(z)$
for $i\in I$.

\item The limits $\lim_{z\to 0}\phi^{i}(z)\one$ for $i\in I$ exist (the limits can be taken
in the topology of $\overline{V}$, but by Lemma \ref{laurent} above,
the existence of these limits mean that the expansions of $\phi^{i}(z)$ have only nonnegative powers).
These elements are either in $V^{0}$ or $V^{1}$. We define $|\phi^{i}|=0$ if
 $\phi^{i}_{-1}\one=\lim_{z\to 0}\phi^{i}(z)\one
\in V^{0}$ and $|\phi^{i}|=1$ if $\phi^{i}_{-1}\one=\lim_{z\to 0}\phi^{i}(z)\one
\in V^{1}$.

\item The vector space $V$ is spanned by elements of the form
$\phi^{i_{1}}_{n_{1}}\cdots \phi^{i_{k}}_{n_{k}}\one$ for $i_{1}, \dots, i_{k}\in I$
and $n_{1}, \dots, n_{k}\in \Z$.

\item  For $v'\in V'$, $v\in V$ and $i_{1}, \dots, i_{k}\in I$, the series
$\langle v', \phi^{i_{1}}(z_{1})\cdots \phi^{i_{k}}(z_{k})v\rangle$
(in fact a Laurent series  in $z_{1}, \dots, z_{k}$ with complex coefficients by Lemma \ref{laurent})
is absolutely convergent in the region $|z_{1}|>\cdots>|z_{k}|>0$ to a
rational function
$R(\langle v', \phi^{i_{1}}(z_{1})\cdots \phi^{i_{k}}(z_{k})v\rangle)$
in $z_{1}, \dots, z_{k}$ with the only possible poles at $z_{i}=0$ for $i=1, \dots, k$ and
$z_{j}=z_{l}$ for $j\ne l$.
In addition, the order of the pole $z_{j}=z_{l}$
is independent of $\phi^{i_{n}}$ for $n\ne j, l$, $v$ and $v'$ and the order of the pole
$z_{j}=0$ is independent of $\phi^{i_{n}}$ for $n\ne j$ and $v'$.

\item  For $v\in V$, $v'\in V'$, $i_{1}, i_{2}\in I$,
$$
R(\langle v', \phi^{i_{1}}(z_{1})\phi^{i_{2}}(z_{2})v\rangle)
=(-1)^{|\phi^{i_{1}}||\phi^{i_{2}}|} R(\langle v', \phi^{i_{2}}(z_{2})
\phi^{i_{1}}(z_{1})v\rangle).
$$

\end{enumerate}

\begin{prop}\label{properties}
The space $V$, the maps $\phi^{i}$ for $i\in I$, $L_{V}(-1)$ and $\one$ have the following
properties:

\begin{enumerate}
\setcounter{enumi}{6}

\item For $a\in \C$ and $i\in I$, $e^{aL_{V}(0)}\phi^{i}(z)e^{-aL_{V}(0)}=e^{a(\swt \phi^{i})}\phi^{i}(e^{a}z)$.

\item $L_{V}(-1) \phi^{i_{1}}_{n_{1}}\cdots \phi^{i_{k}}_{n_{k}}\one={\displaystyle \sum_{j=1}^{k}
\phi^{i_{1}}_{n_{1}}\cdots \phi^{i_{j-1}}_{n_{j-1}}(-n_{j}\phi^{i_{j}}_{n_{j}-1})
\phi^{i_{j+1}}_{n_{j+1}}\cdots \phi^{i_{k}}_{n_{k}}\one}$.

\item For $a\in \C$, $z\in \C^{\times}$ satisfying $|z|>|a|$ and $i\in I$,
$e^{aL_{V}(-1)}\phi^{i}(z)e^{-aL_{V}(-1)}=\phi^{i}(z+a)$ .

\item The operator $L_{V}(-1)$ has weight $1$ and its adjoint $L_{V}(-1)'$ as an operator
on $V'$ has weight $-1$ (the weight of an operator on $V'$ is defined in the same way as
that of an operator on $V$). In particular, $e^{zL_{V}(-1)'}v'\in V'$ for $z\in \C$ and $v'\in V'$.

\item  For $v\in V$, $v'\in V'$ and $\sigma\in S_{k}$,
$$
R(\langle v', \phi^{i_{1}}(z_{1})\cdots \phi^{i_{k}}(z_{k})v\rangle)
=\pm R(\langle v', \phi^{i_{\sigma(1)}}(z_{\sigma(1)})\cdots
\phi^{i_{\sigma(k)}}(z_{\sigma(k)})v\rangle),
$$
where the sign $\pm$ is uniquely determined by $\sigma$
and $|\phi^{i_{1}}|, \dots, |\phi^{i_{k}}|$ (here we omit its
explicit but complicated formula that can be calculated easily for special cases using Condition 6).

\end{enumerate}
\end{prop}
\pf
These properties follow immediately from Conditions 1--6.
\epfv

For $N\in \Z_{+}$, $(z_{i}-z_{j})^{N}$ is a
polynomial in $z_{1}$ and $z_{2}$. We shall use $(z_{i}-z_{j})^{N_{ij}}\phi(z_{i})\phi(z_{j})$
to denote the Laurent series obtained by multiplying the polynomial $(z_{i}-z_{j})^{N_{ij}}$
to the Laurent series $\phi(z_{i})\phi(z_{j})$. We warn the reader that, unless otherwise stated,
$(z_{i}-z_{j})^{N_{ij}}\phi(z_{i})\phi(z_{j})$
is not the Laurent series obtained by multiplying the complex number $(z_{i}-z_{j})^{N_{ij}}$ to the Laurent series $\phi(z_{i})\phi(z_{j})$.

\begin{prop}\label{cond-13}
Let $V=\coprod_{n\in \frac{\Z}{2}}V_{(n)}$ be a $\frac{\Z}{2}$-graded vector space,
$\phi^{i}$ for $i\in I$ analytic maps from
$\C^{\times}$ to $\hom(V, \overline{V})$, $L_{V}(-1)$ an operator on $V$ and
$\one\in V_{(0)}$.
Assume that they satisfy Conditions 1--4. Then Conditions 5 and 6 are equivalent to
the following weak commutativity:
\begin{enumerate}
\setcounter{enumi}{11}
\item  For $i, j\in I$, there exists $N_{ij}\in \Z_{+}$ such that
\begin{equation}\label{locality}
(z_{1}-z_{2})^{N_{ij}}\phi^{i}(z_{1})\phi^{j}(z_{2})=(z_{1}-z_{2})^{N_{ij}}(-1)^{|\phi^{i}||\phi^{j}|}\phi^{j}(z_{2})\phi^{i}(z_{1}).
\end{equation}
\end{enumerate}
In particular, when Conditions 1--4 and Property 12 holds, properties 7--11 also hold.
\end{prop}
\pf
It is clear that Conditions 5 and 6 imply Property 12.

Now we assume that Property 12 holds.
Consider the Laurent series
\begin{equation}\label{properties-1}
\prod_{1\le p<q\le k}(z_{p}-z_{q})^{N_{i_{p}i_{q}}}\langle v', \phi^{i_{1}}(z_{1})\cdots \phi^{i_{k}}(z_{k})v\rangle.
\end{equation}
For $1\le l\le k$, using (\ref{locality}), the Laurent series (\ref{properties-1}) is equal to
\begin{equation}\label{properties-2}
\prod_{1\le p<q\le k}(z_{p}-z_{q})^{N_{i_{p}i_{q}}}\langle v', \phi^{i_{1}}(z_{1})\cdots \phi^{i_{l-1}}(z_{l-1})\phi^{i_{l+1}}(z_{l+1})\cdots \phi^{i_{k}}(z_{k})
\phi^{i_{l}}(z_{l})v\rangle.
\end{equation}
By Lemma \ref{laurent},  (\ref{properties-2})  has only finitely many negative power terms in $z_{l}$.
So the same is true for (\ref{properties-1}).  On the other hand, using (\ref{locality}) again,
(\ref{properties-1}) is equal to
\begin{equation}\label{properties-3}
\prod_{1\le p<q\le k}(z_{p}-z_{q})^{N_{i_{p}i_{q}}}\langle v', \phi^{i_{l}}(z_{l})\phi^{i_{1}}(z_{1})\cdots
\phi^{i_{l-1}}(z_{l-1})\phi^{i_{l+1}}(z_{l+1})\cdots \phi^{i_{k}}(z_{k})
v\rangle.
\end{equation}
By Lemma \ref{laurent} again,  (\ref{properties-3})  has only finitely many positive power terms in $z_{l}$.
So the same is true for (\ref{properties-1}). Thus (\ref{properties-1}) must be a Laurent polynomial in $z_{l}$.
Since this is true for $1\le l\le k$,  (\ref{properties-1}) is a Laurent polynomial in $z_{1}, \dots, z_{k}$.

For fixed $1\le p< q\le k$, the expansion coefficients of
\begin{equation}\label{properties-4}
\langle v', \phi(z_{1})\cdots \phi(z_{k})v\rangle
\end{equation}
as Laurent series in $z_{l}$ for $l\ne p, q$ are of the form
\begin{equation}\label{properties-5}
\langle v', \phi^{i_{1}}_{n_{1}}\cdots \phi^{i_{p-1}}_{n_{p-1}}\phi^{i_{p}}(z_{p})\phi^{i_{p+1}}_{n_{p+1}}\cdots
\phi^{i_{q-1}}_{n_{q-1}}\phi^{i_{q}}(z_{q})\phi^{i_{q+1}}_{n_{q+1}}\cdots \phi^{i_{k}}_{n_{k}}v\rangle
\end{equation}
for $n_{l}\in \Z$, $l\ne p, q$. Clearly (\ref{properties-5}) contains only finitely many
negative powers in $z_{q}$ and finitely many positive powers in $z_{p}$. But we have shown that
when multiplied by $(z_{p}-z_{q})^{N_{pq}}$, it becomes a Laurent polynomial. Thus (\ref{properties-5}) must be
the product of a Laurent polynomial in $z_{p}$ and $z_{q}$ and the expansion of $(z_{p}-z_{q})^{-N_{pq}}$ as a
Laurent series in nonnegative powers of $z_{q}$, or equivalently, in the region $|z_{p}|>|z_{q}|>0$.
Since $p$ and $q$ are arbitrary, we see that
(\ref{properties-4}) is equal to the product of a Laurent polynomial and the expansion of $\prod_{1\le p<q\le k}(z_{p}-z_{q})^{-N_{pq}}$
in the region $|z_{1}|>\cdots >|z_{k}|>0$. This is Condition 5. Condition 6 follows immediately from
Condition 5 in the case $k=2$ and (\ref{locality}).
\epfv

We now define a vertex operator map.
We first give the motivation of this definition. The vertex operator map
we want to define is a map
\begin{eqnarray*}
Y_{V}: \C^{\times}&\to& \hom(V\otimes V, \overline{V}),\nn
z&\mapsto &Y_{V}(\cdot, z)\cdot: u\otimes v\mapsto Y_{V}(u, z)v.
\end{eqnarray*}
We define
$Y_{V}(\phi_{-1}^{i}\one, z)v=\phi^{i}(z)v$ for $i\in I$ and $v\in V$.
The vertex operator map should satisfy the rationality and associativity property.
In particular, we should have
$$R(\langle v', Y_{V}(\phi^{i_{1}}(\xi_{1})\cdots \phi^{i_{k}}(\xi_{k})\one, z)v\rangle)
=R(\langle v', \phi^{i_{1}}(\xi_{1}+z)\cdots \phi^{i_{k}}(\xi_{k}+z)v\rangle)$$
for $i_{1}, \dots, i_{k}\in I$, $v\in V$ and $v'\in V'$.

Motivated by this associativity formula, we define the vertex operator map as follows:
For $v'\in V'$, $v\in V$, $i_{1}, \dots, i_{k}\in I$, $m_{1}, \dots, m_{k}\in \Z$,
we define $Y_{V}$ by
\begin{eqnarray}\label{vo}
\lefteqn{\langle v', Y_{V}(\phi^{i_{1}}_{m_{1}}\cdots \phi^{i_{k}}_{m_{k}}\one, z)
v\rangle}\nn
&&=\res_{\xi_{1}=0}\cdots\res_{\xi_{k}=0}
\xi_{1}^{m_{1}}\cdots\xi_{k}^{m_{k}}
R(\langle v', \phi^{i_{1}}(\xi_{1}+z)\cdots \phi^{i_{k}}(\xi_{k}+z)
v\rangle).
\end{eqnarray}
Note that for a meromorphic function $f(\xi)$, $\res_{\xi=0}f(\xi)$ means expanding $f(\xi)$
as a Laurent series in $0<|\xi|<r$ for $r$ sufficiently small so that no other poles are in this disk
and then taking the coefficient of $\xi^{-1}$. We can also
expand  $f(\xi)$ as a Laurent series in a different region. In general, the coefficient of $\xi^{-1}$
in this Laurent series might be different from $\res_{\xi=0}f(\xi)$.
Also note that the order to take these residues is important. Different orders in general give vertex
operators for different elements.

Since $\overline{V}=(V')^{*}$, for fixed $\phi^{i_{1}}_{m_{1}}\cdots \phi^{i_{k}}_{m_{k}}\one,
v\in V$,
the formula above indeed gives an element
$$Y_{V}(\phi^{i_{1}}_{m_{1}}\cdots \phi^{i_{k}}_{m_{k}}\one, z)
v\in \overline{V}.$$

Since there might be relations among elements of the form
$\phi^{i_{1}}_{m_{1}}\cdots \phi^{i_{k}}_{m_{k}}\one$,
we first have to show that the definition above indeed gives a well-defined map
from $\C^{\times}$ to $\hom(V\otimes V, \overline{V})$. Let $\phi^{0}$ be the map
from $\C^{\times}$ to $\hom(V, \overline{V})$ given by $\phi^{0}(z)=1_{V}$.
Let $\wt \phi^{0}=0$. Then Conditions 1  to 6 and Properties 7 to 12 above still hold for
$\phi^{i}$, $i\in \tilde{I}=I\cup \{0\}$.
Then any relation among
such elements can always be written as
$$\sum_{p=1}^{q}\lambda_{p}\phi^{i^{p}_{1}}_{m^{p}_{1}}\cdots \phi^{i^{p}_{k}}_{m^{p}_{k}}\one=0$$
for some $ i^{p}_{j}\in \tilde{I}$ and $m^{p}_{j}\in \Z$, $p=1, \dots, q$, $j=1, \dots, k$.

\begin{lemma}
If $$\sum_{p=1}^{q}\lambda_{p}\phi^{i^{p}_{1}}_{m^{p}_{1}}\cdots \phi^{i^{p}_{k}}_{m^{p}_{k}}\one=0,$$
then
$$\sum_{p=1}^{q}\lambda_{p}\res_{\xi_{1}=0}\cdots\res_{\xi_{k}=0}
\xi_{1}^{m^{p}_{1}}\cdots\xi_{k}^{m^{p}_{k}}R(\langle v', \phi^{i^{p}_{1}}(\xi_{1}+z)\cdots \phi^{i^{p}_{k}}(\xi_{k}+z)
v\rangle)=0$$
for $v\in V$ and $v'\in V'$.
\end{lemma}
\pf
By Condition 4, we can take $v$ to be of the form $\phi^{j_{1}}_{n_{1}}\cdots \phi^{j_{l}}_{n_{l}}\one$. Moreover,
in this case,
\begin{eqnarray*}
\lefteqn{R(\langle v', \phi^{i^{p}_{1}}(z_{1})\cdots \phi^{i^{p}_{k}}(z_{k})
\phi^{j_{1}}_{n_{1}}\cdots \phi^{j_{l}}_{n_{l}}\one\rangle)}\nn
&&=\res_{\zeta_{1}=0}\cdots\res_{\zeta_{l}=0}
\zeta_{1}^{n_{1}}\cdots\zeta_{k}^{n_{l}}R(\langle v', \phi^{i^{p}_{1}}(z_{1})\cdots \phi^{i^{p}_{k}}(z_{k})
\phi^{j_{1}}(\zeta_{1})\cdots \phi^{j_{l}}(\zeta_{l})\one\rangle).
\end{eqnarray*}
Then
\begin{eqnarray*}
\lefteqn{\res_{\xi_{1}=0}\cdots\res_{\xi_{k}=0}
\xi_{1}^{m^{p}_{1}}\cdots\xi_{k}^{m^{p}_{k}}R(\langle v', \phi^{i^{p}_{1}}(\xi_{1}+z)\cdots \phi^{i^{p}_{k}}(\xi_{k}+z)
v\rangle)}\nn
&&=\res_{\xi_{1}=0}\cdots\res_{\xi_{k}=0}
\xi_{1}^{m^{p}_{1}}\cdots\xi_{k}^{m^{p}_{k}}\res_{\zeta_{1}=0}\cdots\res_{\zeta_{l}=0}
\zeta_{1}^{n_{1}}\cdots\zeta_{l}^{n_{l}}\cdot\nn
&&\quad\quad\quad\quad\quad\cdot R(\langle v', \phi^{i^{p}_{1}}(\xi_{1}+z)\cdots \phi^{i^{p}_{k}}(\xi_{k}+z)
\phi^{j_{1}}(\zeta_{1})\cdots \phi^{j_{l}}(\zeta_{l})\one\rangle)\nn
&&=\prod_{r=1}^{k}\prod_{s=1}^{l}(-1)^{|\phi^{i_{r}}||\phi^{j_{s}}|} \res_{\xi_{1}=0}\cdots\res_{\xi_{k}=0}
\xi_{1}^{m^{p}_{1}}\cdots\xi_{k}^{m^{p}_{k}}\res_{\zeta_{1}=0}\cdots\res_{\zeta_{l}=0}
\zeta_{1}^{n_{1}}\cdots\zeta_{l}^{n_{l}}\cdot\nn
&&\quad\quad\quad\quad\quad\cdot R(\langle v',
\phi^{j_{1}}(\zeta_{1})\cdots \phi^{j_{l}}(\zeta_{l})\phi^{i^{p}_{1}}(\xi_{1}+z)\cdots \phi^{i^{p}_{k}}(\xi_{k}+z)\one\rangle)\nn
&&=\prod_{r=1}^{k}\prod_{s=1}^{l}(-1)^{|\phi^{i_{r}}||\phi^{j_{s}}|}  \res_{\xi_{1}=0}\cdots\res_{\xi_{k}=0}
\xi_{1}^{m^{p}_{1}}\cdots\xi_{k}^{m^{p}_{k}}\res_{\zeta_{1}=0}\cdots\res_{\zeta_{l}=0}
\zeta_{1}^{n_{1}}\cdots\zeta_{l}^{n_{l}}\cdot\nn
&&\quad\quad\quad\quad\quad\cdot R(\langle e^{zL_{V}(-1)'}v',
\phi^{j_{1}}(\zeta_{1}-z)\cdots \phi^{j_{l}}(\zeta_{l}-z)\phi^{i^{p}_{1}}(\xi_{1})\cdots \phi^{i^{p}_{k}}(\xi_{k})\one\rangle)\nn
&&=\prod_{r=1}^{k}\prod_{s=1}^{l}(-1)^{|\phi^{i_{r}}||\phi^{j_{s}}|}  \res_{\zeta_{1}=0}\cdots\res_{\zeta_{l}=0}
\zeta_{1}^{n_{1}}\cdots\zeta_{l}^{n_{l}}\cdot\nn
&&\quad\quad\quad\quad\quad \cdot  R(\langle e^{zL_{V}(-1)'}v',
\phi^{j_{1}}(\zeta_{1}-z)\cdots \phi^{j_{l}}(\zeta_{l}-z)\phi^{i^{p}_{1}}_{m^{p}_{1}}\cdots \phi^{i^{p}_{k}}_{m^{p}_{k}}\one\rangle).
\end{eqnarray*}
Thus
\begin{eqnarray*}
\lefteqn{\sum_{p=1}^{q}\lambda_{p}\res_{\xi_{1}=0}\cdots\res_{\xi_{k}=0}
\xi_{1}^{m^{p}_{1}}\cdots\xi_{k}^{m^{p}_{k}}R(\langle v', \phi^{i^{p}_{1}}(\xi_{1}+z)\cdots \phi^{i^{p}_{k}}(\xi_{k}+z)
v\rangle)}\nn
&&=\sum_{p=1}^{q}\prod_{r=1}^{k}\prod_{s=1}^{l}(-1)^{|\phi^{i_{r}}||\phi^{j_{s}}|}\lambda_{p}\res_{\zeta_{1}=0}\cdots\res_{\zeta_{l}=0}
\zeta_{1}^{n_{1}}\cdots\zeta_{l}^{n_{l}} \cdot\nn
&&\quad\quad\quad\quad\quad \cdot R(\langle e^{zL_{V}(-1)'}v',
\phi^{j_{1}}(\zeta_{1}-z)\cdots \phi^{j_{l}}(\zeta_{l}-z)\phi^{i^{p}_{1}}_{m^{p}_{1}}\cdots \phi^{i^{p}_{k}}_{m^{p}_{k}}\one\rangle)\nn
&&=\prod_{r=1}^{k}\prod_{s=1}^{l}(-1)^{|\phi^{i_{r}}||\phi^{j_{s}}|}\res_{\zeta_{1}=0}\cdots\res_{\zeta_{l}=0}
\zeta_{1}^{n_{1}}\cdots\zeta_{l}^{n_{l}} \cdot\nn
&&\quad\quad\quad\quad\quad \cdot
R\left(\left\langle e^{zL_{V}(-1)'}v',
\phi^{j_{1}}(\zeta_{1}-z)\cdots \phi^{j_{l}}(\zeta_{l}-z)\left(\sum_{p=1}^{q}\lambda_{p}\phi^{i^{p}_{1}}_{m^{p}_{1}}\cdots
\phi^{i^{p}_{k}}_{m^{p}_{k}}\one\right)\right\rangle\right)\nn
&&=0,
\end{eqnarray*}
proving the lemma.
\epfv

From this lemma, we see that the vertex operator map $Y_{V}$ is well defined.
We are now ready to formulate and prove the main result of this section.

\begin{thm}\label{first-const}
Let $V=\coprod_{n\in \frac{\Z}{2}}V_{(n)}$ be a $\frac{\Z}{2}$-graded vector space,
$\phi^{i}$ for $i\in I$ maps from
$\C^{\times}$ to $\hom(V, \overline{V})$, $L_{V}(-1)$ an operator on $V$ and
$\one\in V_{(0)}$. Assume that they satisfy Conditions 1--6.
Then the triple $(V, Y_{V}, \one)$ is a grading-restricted vertex algebra generated
by $\phi^{i}_{-1}\one$ for $i\in I$.
Moreover, this is the unique grading-restricted vertex algebra structure on $V$ with the vacuum $\one$
such that $Y(\phi^{1}_{-1}\one, z)=\phi^{i}(z)$
for $i\in I$.
\end{thm}
\pf
The vertex operator map $Y_{V}$ is clearly analytic.
The grading-restriction axiom is by assumption satisfied. The $L(-1)$-bracket formula follows from
Condition 1 and the definition of $Y_{V}$. The identity property and the creation property
also follow from of the definition of $Y_{V}$.

Let $L_{V}(0)'$ be the adjoint operator of $L_{V}(0)$.
For $v'\in V'$, $v\in V$, $ i_{1}, \dots, i_{k}\in I$ and $n_{1}, \dots, n_{k}\in \Z$, $a\in \C^{\times}$
\begin{eqnarray*}
\lefteqn{\langle  v', a^{L_{V}(0)}Y_{V}(\phi^{i_{1}}_{n_{1}}\cdots \phi^{i_{k}}_{n_{k}}\one, z)a^{-L_{V}(0)}v\rangle}\nn
&&=\langle  a^{L_{V}(0)'}v', Y_{V}(\phi^{i_{1}}_{n_{1}}\cdots \phi^{i_{k}}_{n_{k}}\one, z)a^{-L_{V}(0)}v\rangle\nn
&&=\res_{\xi_{1}=0}\cdots\res_{\xi_{k}=0}\xi_{1}^{n_{1}}\cdots \xi_{k}^{n_{k}}
R(\langle  a^{L_{V}(0)'}v', \phi^{i_{1}}(\xi_{1}+z)\cdots \phi^{i_{k}}(\xi_{k}+z)a^{-L_{V}(0)}v\rangle)\nn
&&=\res_{\xi_{1}=0}\cdots\res_{\xi_{k}=0}\xi_{1}^{n_{1}}\cdots \xi_{k}^{n_{k}}
R(\langle  v', a^{L_{V}(0)}\phi^{i_{1}}(\xi_{1}+z)\cdots \phi^{i_{k}}(\xi_{k}+z)a^{-L_{V}(0)}v\rangle)\nn
&&=\res_{\xi_{1}=0}\cdots\res_{\xi_{k}=0}\xi_{1}^{n_{1}}\cdots \xi_{k}^{n_{k}}a^{\swt \phi^{i_{1}}+\cdots
\swt \phi^{i_{k}}}
R(\langle  v', \phi^{i_{1}}(a\xi_{1}+az)\cdots \phi^{i_{k}}(a\xi_{k}+az)v\rangle)\nn
&&=\res_{\zeta_{1}=0}\cdots\res_{\zeta_{k}=0}\zeta_{1}^{n_{1}}\cdots \zeta_{k}^{n_{k}}a^{\swt \phi^{i_{1}}+\cdots
\swt \phi^{i_{k}}-k-n_{1}-\cdots -n_{k}}\cdot\nn
&&\quad\quad\quad\quad\quad\quad\quad\quad\quad\quad\quad\quad\quad\quad\cdot
R(\langle  v', \phi^{i_{1}}(\zeta_{1}+az)\cdots \phi^{i_{k}}(\zeta_{k}+az)v\rangle)\nn
&&=\langle  v', Y_{V}(a^{L_{V}(0)}\phi^{i_{1}}_{n_{1}}\cdots \phi^{i_{k}}_{n_{k}}\one , az)v\rangle).
\end{eqnarray*}
This formula implies the $L(0)$-bracket formula.

From Condition 2 and the definition of $Y_{V}$, we obtain
$$\frac{d}{dz}Y_{V}(\phi^{i_{1}}_{n_{1}}\cdots \phi^{i_{k}}_{n_{k}}\one, z)
=[L_{V}(-1), Y_{V}(\phi^{i_{1}}_{n_{1}}\cdots \phi^{i_{k}}_{n_{k}}\one, z)].$$
From Property 8 and the definition of $Y_{V}$, we obtain
$$\frac{d}{dz}Y_{V}(\phi^{i_{1}}_{n_{1}}\cdots \phi^{i_{k}}_{n_{k}}\one, z)
=Y_{V}(L_{V}(-1)\phi^{i_{1}}_{n_{1}}\cdots \phi^{i_{k}}_{n_{k}}\one, z).$$
Applying both sides of this formula to $\one$, taking the limit $z\to 0$
and then using the creation property, we obtain
$$L_{V}(-1)\phi^{i_{1}}_{n_{1}}\cdots \phi^{i_{k}}_{n_{k}}\one=\lim_{z\to 0}
\frac{d}{dz}Y_{V}(\phi^{i_{1}}_{n_{1}}\cdots \phi^{i_{k}}_{n_{k}}\one, z)\one.$$
The $L(-1)$-derivative property is proved.

Let $\{e_{n}\}_{n\in \Z}$ be a homogeneous basis of $V$ and $\{e_{n}'\}_{n\in \Z}$ its dual basis in $V'$.
Then we have
\begin{eqnarray}\label{prod}
\lefteqn{\langle v', Y_{V}(\phi^{i_{1}}_{n_{1}}\cdots \phi^{i_{k}}_{n_{k}}\one, z_{1})
Y_{V}(\phi^{j_{1}}_{m_{1}}\cdots \phi^{j_{l}}_{m_{l}}\one, z_{2})v\rangle}\nn
&&=\sum_{n\in \Z}\langle v', Y_{V}(\phi^{i_{1}}_{n_{1}}\cdots \phi^{i_{k}}_{n_{k}}\one, z_{1})e_{n}\rangle
\langle e'_{n}, Y_{V}(\phi^{j_{1}}_{m_{1}}\cdots \phi^{j_{l}}_{m_{l}}\one, z_{2})v\rangle\nn
&&=\sum_{n\in \Z}\res_{\zeta_{1}=0}\cdots\res_{\zeta_{k}=0}\zeta_{1}^{n_{1}}\cdots \zeta_{k}^{n_{k}}
\res_{\xi_{1}=0}\cdots\res_{\xi_{l}=0}\xi_{1}^{m_{1}}\cdots \xi_{l}^{m_{l}}\cdot\nn
&&\quad\quad\quad\quad\cdot R(\langle v', \phi^{i_{1}}(\zeta_{1}+z_{1})\cdots
 \phi^{i_{k}}(\zeta_{k}+z_{1})e_{n}\rangle)
R(\langle e'_{n}, \phi^{j_{1}}(\xi_{1}+z_{2})\cdots \phi^{j_{l}}(\xi_{l}+z_{2})v\rangle)\nn
&&=\res_{\zeta_{1}=0}\cdots\res_{\zeta_{k}=0}\zeta_{1}^{n_{1}}\cdots \zeta_{k}^{n_{k}}
\res_{\xi_{1}=0}\cdots\res_{\xi_{l}=0}\xi_{1}^{m_{1}}\cdots \xi_{l}^{m_{l}}\cdot\nn
&&\quad\quad\quad\quad\cdot \sum_{n\in \Z}R(\langle v', \phi^{i_{1}}(\zeta_{1}+z_{1})\cdots
 \phi^{i_{k}}(\zeta_{k}+z_{1})e_{n}\rangle)
R(\langle e'_{n}, \phi^{j_{1}}(\xi_{1}+z_{2})\cdots \phi^{j_{l}}(\xi_{l}+z_{2})v\rangle).\nn
\end{eqnarray}

By Condition 5, when $|z_{1}|>\cdots >|z_{k+l}|>0$,
\begin{eqnarray}\label{prod-1}
\lefteqn{\sum_{n\in \Z}R(\langle v', \phi^{i_{1}}(z_{1})\cdots
 \phi^{i_{k}}(z_{k})e_{n}\rangle)
R(\langle e'_{n}, \phi^{j_{1}}(z_{k+1})\cdots \phi^{j_{l}}(z_{k+l})v\rangle)}\nn
&&=\sum_{n\in \Z}\langle v', \phi^{i_{1}}(z_{1})\cdots
 \phi^{i_{k}}(z_{k})e_{n}\rangle
\langle e'_{n}, \phi^{j_{1}}(z_{k+1})\cdots \phi^{j_{l}}(z_{k+l})v\rangle\nn
&&=\langle v', \phi^{i_{1}}(z_{1})\cdots
 \phi^{i_{k}}(z_{k})\phi^{j_{1}}(z_{k+1})\cdots \phi^{j_{l}}(z_{k+l})v\rangle
\end{eqnarray}
is absolutely convergent to the rational function
\begin{equation}\label{prod-2}
R(\langle v', \phi^{i_{1}}(z_{1})\cdots
 \phi^{i_{k}}(z_{k})\phi^{j_{1}}(z_{k+1})\cdots \phi^{j_{l}}(z_{k+l})v\rangle)
\end{equation}
in $z_{1}, \dots, z_{k+l}$. On the other hand, since the only possible poles of (\ref{prod-2})
are $z_{i}-z_{j}=0$ for $i\ne j$ and $z_{i}=0$, there is a unique expansion of such a rational function
in the region $|z_{1}|, \dots, |z_{k}|>|z_{k+1}|, \dots, |z_{k+l}|>0$, $z_{i}\ne z_{j}$ for $i\ne j$,
$i, j=1, \dots, k$ and $i, j=k+1, \dots, k+l$ such that each term
is a product of two rational functions, one in $z_{1}, \dots, z_{k}$ and the other in $z_{k+1}, \dots,
z_{k+l}$. Since the left-hand side of (\ref{prod-1}) is a series of the same form
and is absolutely convergent
in the region $|z_{1}|>\cdots >|z_{k+l}|>0$ to  (\ref{prod-2}), it must be absolutely convergent
in the larger region $|z_{1}|, \dots,|z_{k}|>|z_{k+1}|, \dots, |z_{k+l}|>0$,
$z_{i}\ne z_{j}$ for $i\ne j$,
$i, j=1, \dots, k$ and $i, j=k+1, \dots, k+l$ to (\ref{prod-2}).

Substituting $\zeta_{i}+z_{1}$ for $z_{i}$ for $i=1, \dots, k$ and $\xi_{j}+z_{2}$ for $z_{k+j}$
for $j=1, \dots, l$, we see that
$$ \sum_{n\in \Z}R(\langle v', \phi^{i_{1}}(\zeta_{1}+z_{1})\cdots
 \phi^{i_{k}}(\zeta_{k}+z_{1})e_{n}\rangle)
R(\langle e'_{n}, \phi^{j_{1}}(\xi_{1}+z_{2})\cdots \phi^{j_{l}}(\xi_{l}+z_{2})v\rangle)$$
is absolutely convergent to
$$R(\langle v', \phi^{i_{1}}(\zeta_{1}+z_{1})\cdots \phi^{i_{k}}(\zeta_{k}+z_{1})
\phi^{j_{1}}(\xi_{1}+z_{2})\cdots \phi^{j_{l}}(\xi_{l}+z_{2})v\rangle)$$
when $|\zeta_{1}+z_{1}|, \dots, |\zeta_{k}+z_{1}|>|\xi_{1}+z_{2}|, \dots, |\xi_{l}+z_{2}|>0$,
$\zeta_{i}\ne \zeta_{j}$  for $i, j=1, \dots, k$ and $\xi_{i}\ne \xi_{j}$ for $i, j=1, \dots, l$. When
$|z_{1}|>|z_{2}|>0$, we can always find sufficiently small neighborhood of $0$ such that when
$\zeta_{1}, \dots, \zeta_{k}, \xi_{1}, \dots, \xi_{l}$ are in this neighborhood,
$|\zeta_{1}+z_{1}|, \dots, |\zeta_{k}+z_{1}|>|\xi_{1}+z_{2}|, \dots, |\xi_{l}+z_{2}|>0$
holds. Thus we see that when $|z_{1}|>|z_{2}|>0$, the right-hand side of (\ref{prod}) is absolutely convergent to
\begin{eqnarray}\label{prod-3}
\lefteqn{\res_{\zeta_{1}=0}\cdots\res_{\zeta_{k}=0}\zeta_{1}^{n_{1}}\cdots \zeta_{k}^{n_{k}}
\res_{\xi_{1}=0}\cdots\res_{\xi_{l}=0}\xi_{1}^{m_{1}}\cdots \xi_{l}^{m_{l}}\cdot}\nn
&&\quad\quad\quad\quad\quad\cdot R(\langle v', \phi^{i_{1}}(\zeta_{1}+z_{1})\cdots \phi^{i_{k}}(\zeta_{k}+z_{1})
\phi^{j_{1}}(\xi_{1}+z_{2})\cdots \phi^{j_{l}}(\xi_{l}+z_{2})v\rangle).
\end{eqnarray}
This is a rational function in $z_{1}$ and $z_{2}$ with the only
possible poles at $z_{1}, z_{2}=0$ and $z_{1}=z_{2}$. In particular, the left-hand side of
(\ref{prod}), that is,
\begin{equation}\label{prod-4}
\langle v', Y_{V}(\phi^{i_{1}}_{n_{1}}\cdots \phi^{i_{k}}_{n_{k}}\one, z_{1})
Y_{V}(\phi^{j_{1}}_{m_{1}}\cdots \phi^{j_{l}}_{m_{l}}\one, z_{2})v\rangle,
\end{equation}
is absolutely convergent in the region $|z_{1}|>|z_{2}|>0$ to this rational
function.

We have proved the rationality of the product of two vertex operators. We are ready to
prove the commutativity.
The calculation above also shows that
\begin{equation}\label{prod-5}
\langle v',
Y_{V}(\phi^{j_{1}}_{m_{1}}\cdots \phi^{j_{l}}_{m_{l}}\one, z_{2})
Y_{V}(\phi^{i_{1}}_{n_{1}}\cdots \phi^{i_{k}}_{n_{k}}\one, z_{1})v\rangle
\end{equation}
is absolutely convergent to the rational function
\begin{eqnarray}\label{prod-6}
\lefteqn{\res_{\xi_{1}=0}\cdots\res_{\xi_{l}=0}\xi_{1}^{m_{1}}\cdots \xi_{l}^{m_{l}}
\res_{\zeta_{1}=0}\cdots\res_{\zeta_{k}=0}\zeta_{1}^{n_{1}}\cdots \zeta_{k}^{n_{k}}\cdot}\nn
&&\quad\quad\quad\quad\quad\quad\cdot R(\langle v',
\phi^{j_{1}}(\xi_{1}+z_{2})\cdots \phi^{j_{l}}(\xi_{l}+z_{2})
\phi^{i_{1}}(\zeta_{1}+z_{1})\cdots \phi^{i_{k}}(\zeta_{k}+z_{1})v\rangle),\nn
\end{eqnarray}
in the regions $|z_{2}|>|z_{1}|>0$, respectively. By Property 11, the rational functions
(\ref{prod-3}) and (\ref{prod-6}) multiplied by
$$\prod_{r=1}^{k}\prod_{s=1}^{l}(-1)^{|\phi^{i_{r}}||\phi^{j_{s}}|}=(-1)^{|\phi^{i_{1}}_{n_{1}}\cdots \phi^{i_{k}}_{n_{k}}\one|
|\phi^{j_{1}}_{m_{1}}\cdots \phi^{j_{l}}_{m_{l}}\one|}$$
are equal.
Thus (\ref{prod-4}) and (\ref{prod-5}) multiplied by the sign
$(-1)^{|\phi^{i_{1}}_{n_{1}}\cdots \phi^{i_{k}}_{n_{k}}\one|
|\phi^{j_{1}}_{m_{1}}\cdots \phi^{j_{l}}_{m_{l}}\one|}$
are
absolutely convergent in the regions $|z_{1}|>|z_{2}|>0$ and $|z_{2}|>|z_{1}|>0$,
respectively, to a common rational function with the only possible poles at $z_{1}=z_{2}$,
$z_{1}=0$ and $z_{2}=0$.

We now prove the associativity. For $i_{1}, \dots, i_{k}, j_{1}, \dots, j_{l}\in I$, $m_{1}, \dots, m_{l}\in \Z$,
$v\in V$ and $v'\in V'$,
using the expansion of $\phi^{i_{1}}(\xi_{1}), \dots, \phi^{i_{k}}(\xi_{k})$ and the
definition of $Y_{V}$, we have
\begin{eqnarray}\label{phi-assoc}
\lefteqn{\langle v', Y_{V}(\phi^{i_{1}}(z_{1})\cdots \phi^{i_{k}}(z_{k})
\phi^{j_{1}}_{m_{1}}\cdots \phi^{j_{l}}_{m_{l}}\one, z)v\rangle}\nn
&&=\sum_{p_{1}, \dots, p_{k}\in \Z}\langle v', Y_{V}(\phi^{i_{1}}_{p_{1}}\cdots \phi^{i_{k}}_{p_{k}}
\phi^{j_{1}}_{m_{1}}\cdots \phi^{j_{l}}_{m_{l}}\one, z)v\rangle
z_{1}^{-p_{1}-1}\cdots z_{k}^{-p_{k}-1}\nn
&&=\sum_{p_{1}, \dots, p_{k}\in \Z}\res_{\zeta_{1}=0}\cdots\res_{\zeta_{k}=0}\zeta_{1}^{p_{1}}\cdots \zeta_{k}^{p_{k}}
\res_{\xi_{1}=0}\cdots\res_{\xi_{l}=0}\xi_{1}^{m_{1}}\cdots \xi_{l}^{m_{l}}\cdot \nn
&&\quad\quad\quad\quad\cdot
R(\langle v', \phi^{i_{1}}(\zeta_{1}+z)\cdots \phi^{i_{k}}(\zeta_{k}+z)\phi^{j_{1}}(\xi_{1}+z)\cdots \phi^{j_{l}}(\xi_{l}+z)v\rangle)
z_{1}^{-p_{1}-1}\cdots z_{k}^{-p_{k}-1}.\nn
\end{eqnarray}
We now expand
$$R(\langle v', \phi^{i_{1}}(\zeta_{1}+z)\cdots \phi^{i_{k}}(\zeta_{k}+z)\phi^{j_{1}}(\xi_{1}+z)\cdots \phi^{j_{l}}(\xi_{l}+z)v\rangle)$$
 as a Laurent series
$\sum_{l\in \Z}f_{l}(\zeta_{1}, \dots, \zeta_{k-1}, \xi_{1}, \dots, \xi_{l}, z) \zeta_{k}^{-l-1}$ in $\zeta_{k}$
in the region $|z|, |\zeta_{1}|, \dots, |\zeta_{k-1}|>|\zeta_{k}|>|\xi_{1}|, \dots, |\xi_{l}|$, where
$f_{l}(\zeta_{1}, \dots, \zeta_{k-1},  \xi_{1}, \dots, \xi_{l}, z)$ are rational functions in $\zeta_{1}, \dots, \zeta_{k-1},
 \xi_{1}, \dots, \xi_{l}$ and $z$.
Then in the region that the Laurent series expansion holds, we have
\begin{eqnarray}\label{phi-assoc-1}
\lefteqn{\sum_{p_{k}\in \Z}\res_{\zeta_{k}=0}\zeta_{k}^{p_{k}}
\left(\sum_{l\in \Z}f_{l}(\zeta_{1}, \dots, \zeta_{k-1}, \xi_{1}, \dots, \xi_{l}, z) \zeta_{k}^{-l-1}\right)
 z_{k}^{-p_{k}-1}}\nn
&&=\sum_{p_{k}\in \Z}f_{p_{k}}(\zeta_{1}, \dots, \zeta_{k-1}, \xi_{1}, \dots, \xi_{l}, z)  z_{k}^{-p_{k}-1}\nn
&&=R(\langle v', \phi^{i_{1}}(\zeta_{1}+z)\cdots \phi^{i_{k-1}}(\zeta_{k-1}+z) \phi^{i_{k}}(z_{k}+z)
\phi^{j_{1}}(\xi_{1}+z)\cdots \phi^{j_{l}}(\xi_{l}+z)v\rangle).\nn
\end{eqnarray}
Repeating this step for the variables $\zeta_{k-1}, \dots, \zeta_{1}$, we see that the right-hand side of
(\ref{phi-assoc}) is equal to the expansion of
\begin{equation}\label{phi-assoc-2}
\res_{\xi_{1}=0}\cdots\res_{\xi_{l}=0}\xi_{1}^{m_{1}}\cdots \xi_{l}^{m_{l}}
R(\langle v', \phi^{i_{1}}(z_{1}+z)\cdots \phi^{i_{k}}(z_{k}+z)\phi^{j_{1}}(\xi_{1}+z)\cdots \phi^{j_{l}}(\xi_{l}+z)v\rangle)
\end{equation}
as a  Laurent series in $z_{1}\dots, z_{k}$ in the region $|z|>|z_{1}|>\cdots>|z_{k}|>0$.
Thus the left-hand side of (\ref{phi-assoc})
is absolutely convergent to (\ref{phi-assoc-2}) in the region for this Laurent series expansion.
In particular,  in the region $|z|>|z_{1}|>\cdots>|z_{k}|>0$,
\begin{eqnarray}\label{phi-assoc-3}
\lefteqn{\langle v', Y_{V}(\phi^{i_{1}}(z_{1})\cdots \phi^{i_{k}}(z_{k})
\phi^{j_{1}}_{m_{1}}\cdots \phi^{j_{l}}_{m_{l}}\one, z)v\rangle}\nn
&&=\res_{\xi_{1}=0}\cdots\res_{\xi_{l}=0}\xi_{1}^{m_{1}}\cdots \xi_{l}^{m_{l}}\cdot\nn
&&\quad\quad\quad\quad\quad\quad\cdot
R(\langle v', \phi^{i_{1}}(z_{1}+z)\cdots \phi^{i_{k}}(z_{k}+z)\phi^{j_{1}}(\xi_{1}+z)\cdots \phi^{j_{l}}(\xi_{l}+z)v\rangle).
\end{eqnarray}

Now we have
\begin{eqnarray}\label{iter}
\lefteqn{\langle v', Y_{V}(Y_{V}(\phi^{i_{1}}_{n_{1}}\cdots \phi^{i_{k}}_{n_{k}}\one, z_{1}-z_{2})
\phi^{j_{1}}_{m_{1}}\cdots \phi^{j_{l}}_{m_{l}}\one, z_{2})v\rangle}\nn
&&=\sum_{n\in \Z}\langle v', Y_{V}(e_{n}, z_{2})v\rangle\langle e_{n}',
Y_{V}(\phi^{i_{1}}_{n_{1}}\cdots \phi^{i_{k}}_{n_{k}}\one, z_{1}-z_{2})
\phi^{j_{1}}_{m_{1}}\cdots \phi^{j_{l}}_{m_{l}}\one\rangle\nn
&&=\sum_{n\in \Z}\langle v', Y_{V}(e_{n}, z_{2})v\rangle
\res_{\zeta_{1}=0}\cdots\res_{\zeta_{k}=0}\zeta_{1}^{n_{1}}\cdots \zeta_{k}^{n_{k}}\cdot\nn
&&\quad\quad\quad\quad\quad\quad\cdot R(\langle e_{n}',\phi^{i_{1}}(\zeta_{1}+z_{1}-z_{2})\cdots \phi^{i_{k}}(\zeta_{k}+z_{1}-z_{2})
\phi^{j_{1}}_{m_{1}}\cdots \phi^{j_{l}}_{m_{l}}\one\rangle).
\end{eqnarray}
But by (\ref{phi-assoc-3}), in the region $|z_{2}|>|\zeta_{1}+z_{1}-z_{2}|>\cdots>|\zeta_{k}+z_{1}-z_{2}|>0$, we have
\begin{eqnarray}\label{iter-1}
\lefteqn{\sum_{n\in \Z}\langle v', Y_{V}(e_{n}, z_{2})v\rangle
\langle e_{n}',\phi^{i_{1}}(\zeta_{1}+z_{1}-z_{2})\cdots \phi^{i_{k}}(\zeta_{k}+z_{1}-z_{2})
\phi^{j_{1}}_{m_{1}}\cdots \phi^{j_{l}}_{m_{l}}\one\rangle}\nn
&&=\langle v', Y_{V}(\phi^{i_{1}}(\zeta_{1}+z_{1}-z_{2})\cdots \phi^{i_{k}}(\zeta_{k}+z_{1}-z_{2})
\phi^{j_{1}}_{m_{1}}\cdots \phi^{j_{l}}_{m_{l}}\one, z_{2})v\rangle\nn
&&=\res_{\xi_{1}=0}\cdots\res_{\xi_{l}=0}\xi_{1}^{m_{1}}\cdots \xi_{l}^{m_{l}}\cdot\nn
&&\quad\quad\quad\quad\quad\quad\cdot
R(\langle v', \phi^{i_{1}}(\zeta_{1}+z_{1})\cdots \phi^{i_{k}}(\zeta_{k}+z_{1})\phi^{j_{1}}(\xi_{1}+z_{2})\cdots \phi^{j_{l}}(\xi_{l}+z_{2})v\rangle).\nn
\end{eqnarray}
The right-hand side of (\ref{iter-1}) is a rational function in $\zeta_{1},
\dots, \zeta_{k}$, $z_{1}$ and $z_{2}$ with the only possible poles
$\zeta_{i}-\zeta_{j}=0$, for $i\ne j$, $\zeta_{i}+z_{1}=0$, $\zeta_{i}+z_{1}-z_{2}=0$ and $z_{2}=0$. There is a unique expansion of such a rational function
in the region $|z_{2}|>|\zeta_{1}+z_{1}-z_{2}|, \dots, |\zeta_{k}+z_{1}-z_{2}|>0$,
$\zeta_{i}\ne \zeta_{j}$ for $i\ne j$,
$i, j=1, \dots, k$, such that each term
is a product of two rational functions, one in $z_{2}$ and the other in $\zeta_{1},
\dots, \zeta_{k}$ and $z_{1}$. Since
$$\sum_{n\in \Z}\langle v', Y_{V}(e_{n}, z_{2})v\rangle
R(\langle e_{n}',\phi^{i_{1}}(\zeta_{1}+z_{1}-z_{2})\cdots \phi^{i_{k}}(\zeta_{k}+z_{1}-z_{2})
\phi^{j_{1}}_{m_{1}}\cdots \phi^{j_{l}}_{m_{l}}\one\rangle)$$
is a series of the same form  and is equal to
the left-hand side of (\ref{iter-1})  in the region
$|z_{2}|>|\zeta_{1}+z_{1}-z_{2}|>\cdots>|\zeta_{k}+z_{1}-z_{2}|>0$,
it must be absolutely convergent to the right-hand side of (\ref{iter-1}) in the larger region
$|z_{2}|>|\zeta_{1}+z_{1}-z_{2}|, \dots, |\zeta_{k}+z_{1}-z_{2}|>0$. Thus we obtain
\begin{eqnarray}\label{iter-2}
\lefteqn{\sum_{n\in \Z}\langle v', Y_{V}(e_{n}, z_{2})v\rangle
R(\langle e_{n}',\phi^{i_{1}}(\zeta_{1}+z_{1}-z_{2})\cdots \phi^{i_{k}}(\zeta_{k}+z_{1}-z_{2})
\phi^{j_{1}}_{m_{1}}\cdots \phi^{j_{l}}_{m_{l}}\one\rangle)}\nn
&&=\res_{\xi_{1}=0}\cdots\res_{\xi_{l}=0}\xi_{1}^{m_{1}}\cdots \xi_{l}^{m_{l}}\cdot\nn
&&\quad\quad\quad\quad\quad\quad\cdot
R(\langle v', \phi^{i_{1}}(\zeta_{1}+z_{1})\cdots \phi^{i_{k}}(\zeta_{k}+z_{1})\phi^{j_{1}}(\xi_{1}+z_{2})\cdots \phi^{j_{l}}(\xi_{l}+z_{2})v\rangle)\nn
\end{eqnarray}
in the region $|z_{2}|>|\zeta_{1}+z_{1}-z_{2}|, \dots, |\zeta_{k}+z_{1}-z_{2}|>0$.
Thus when $|z_{2}|>|z_{1}-z_{2}|>0$, the right-hand side of (\ref{iter}) is absolutely convergent to
\begin{eqnarray}
\lefteqn{\res_{\zeta_{1}=0}\cdots\res_{\zeta_{k}=0}\zeta_{1}^{n_{1}}\cdots \zeta_{k}^{n_{k}}
\res_{\xi_{1}=0}\cdots\res_{\xi_{l}=0}\xi_{1}^{m_{1}}\cdots \xi_{l}^{m_{l}}\cdot}\nn
&&\quad\quad\quad\quad\quad\quad\cdot R(\langle v', \phi^{i_{1}}(\zeta_{1}+z_{1})\cdots \phi^{i_{k}}(\zeta_{k}+z_{1})
\phi^{j_{1}}(\xi_{1}+z_{2})\cdots \phi^{j_{l}}(\xi_{l}+z_{2})v\rangle),\nn
\end{eqnarray}
which is proved above to be equal to the left hand side of (\ref{prod}) in the region $|z_{1}|>|z_{2}|>0$. The associativity is proved.

To prove the uniqueness, we need only show that any grading-restricted vertex superalgebra
structure on $V$ with the vacuum $\one$
must have the vertex operator map defined by (\ref{vo}). But this is clear from the motivation that we discussed before
the definition (\ref{vo}) of the vertex operator map $Y_{V}$.
\epfv

We call the grading-restricted vertex algebra given in Theorem \ref{first-const}
the {\it grading-restricted vertex algebra generated by $\phi^{i}$, $i\in I$}. The maps
$\phi^{i}$, $i\in I$, are called the {\it generating fields} of the grading-restricted vertex algebra $V$.

A consequence follows immediately from Proposition \ref{cond-13} and Theorem \ref{first-const}
is the following result:

\begin{cor}\label{mp}
Let $V=\coprod_{n\in \frac{\Z}{2}}V_{(n)}$ be a $\frac{\Z}{2}$-graded vector space,
$\phi^{i}$ for $i\in I$ maps from
$\C^{\times}$ to $\hom(V, \overline{V})$, $L_{V}(-1)$ an operator on $V$ and
$\one\in V_{(0)}$. Assume that they satisfy Conditions 1--4 and Property 12.
Then the triple $(V, Y_{V}, \one)$ is a grading-restricted vertex algebra generated
by $\phi^{i}_{-1}\one$ for $i\in I$.
Moreover, this is the unique grading-restricted vertex algebra structure on $V$ with the vacuum $\one$
such that $Y(\phi^{i}_{-1}\one, z)=\phi^{i}(z)$
for $i\in I$. \epf
\end{cor}

\begin{rema}
{\rm From the uniqueness part of Corollary \ref{mp}, we see that the class of grading-restricted
vertex algebras given in Corollary \ref{mp} is the same as the class given by Meurman and Primc in
\cite{MP} using a different construction.
In particular,
the two expressions for the vertex operator maps in the present paper
and in \cite{MP}  must be equal.   But the constructions are completely different.
First, the definitions of vertex operator maps are different. As in all the other constructions
of the second type (see the introduction
of the present paper), the vertex operator map in \cite{MP}
is defined using the residue of  $x_{1}$ in the Jacobi identity for vertex operator algebras.
Our construction uses a completely different definition motivated by the associativity.
Second, the proofs are necessarily different. As in most of the other constructions
of the second type, the proof in \cite{MP}  is based on a result
stating that if a generating field is local with two other fields (or satisfies certain
identities together with two other fields), then this field is also local
with the fields generated by the other two fields (or satisfies certain
identities together with the fields generated by the other two fields). Our proof does not use and does not need to
use such a result.}
\end{rema}

\renewcommand{\theequation}{\thesection.\arabic{equation}}
\renewcommand{\thethm}{\thesection.\arabic{thm}}
\setcounter{equation}{0}
\setcounter{thm}{0}
\section{The second construction}

In this section, we give our second construction.
It gives us a quasi-vertex operator
algebra structure on the direct sum of a quasi-vertex operator algebra and
a module satisfying suitable conditions.

Let $(V, Y_{V}, \one, L_{V}(1))$ be a simple quasi-vertex operator algebra.
Assume that $V$-modules are all completely reducible and $\R$-graded. Also assume  that intertwining
operators among $V$-modules satisfy the associativity property, that is,
for $V$-modules $W_{1}, \dots, W_{5}$ and  intertwining operators $\mathcal{Y}_{1}$ and $\mathcal{Y}_{2}$ of types
${W_{4}\choose W_{1}W_{5}}$ and ${W_{5}\choose W_{2}W_{3}}$, respectively,
there exist a $V$-module $W_{6}$ and intertwining operators $\mathcal{Y}_{3}$ and
$\mathcal{Y}_{4}$ of types ${W_{4}\choose W_{6}W_{3}}$ and ${W_{6}\choose W_{1}W_{2}}$,
respectively, such that for $w_{4}'\in W_{4}'$, $w_{1}\in W_{1}$, $w_{2}\in W_{2}$ and $w_{3}\in W_{3}$,
\begin{equation}\label{assoc-int}
\langle w_{4}', \mathcal{Y}_{1}(w_{1}, z_{1})\mathcal{Y}_{2}(w_{2}, z_{2})w_{3}\rangle
=\langle w_{4}', \mathcal{Y}_{3}(\mathcal{Y}_{4}(w_{1}, z_{1}-z_{2})w_{2}, z_{2})w_{3}\rangle
\end{equation}
when $|z_{1}|>|z_{2}|>|z_{1}-z_{2}|>0$.
See \cite{H4} for conditions on $V$ such that
this associativity holds.

For simplicity, when $|z_{1}|>|z_{2}|>0$, we use
$\mathcal{Y}_{1}(w_{1}, z_{1})\mathcal{Y}_{2}(w_{2}, z_{2})$ to also denote the
element of $\hom(W_{3}, \overline{W}_{4})$ obtained from the sum of the series
$\langle w_{4}', \mathcal{Y}_{1}(w_{1}, z_{1})\mathcal{Y}_{2}(w_{2}, z_{2})w_{3}\rangle$ for
$w_{4}'\in W_{4}'$ and $w_{3}\in W_{3}$.  Similarly, when $|z_{2}|>|z_{1}-z_{2}|>0$, we have an element
$\mathcal{Y}_{3}(\mathcal{Y}_{4}(w_{1}, z_{1}-z_{2})w_{2}, z_{2})\in \hom(W_{3}, \overline{W}_{4})$.
Then when $|z_{1}|>|z_{2}|>|z_{1}-z_{2}|>0$, (\ref{assoc-int})   for all
$w_{4}'\in W_{4}'$ and $w_{3}\in W_{3}$ can be written as
$$\mathcal{Y}_{1}(w_{1}, z_{1})\mathcal{Y}_{2}(w_{2}, z_{2})
=\mathcal{Y}_{3}(\mathcal{Y}_{4}(w_{1}, z_{1}-z_{2})w_{2}, z_{2}).$$
We shall also say that
$$\mathcal{Y}_{1}(w_{1}, z_{1})\mathcal{Y}_{2}(w_{2}, z_{2})$$
defined on $|z_{1}|>|z_{2}|>0$ and
$$\mathcal{Y}_{3}(\mathcal{Y}_{4}(w_{1}, z_{1}-z_{2})w_{2}, z_{2})$$
defined on $|z_{2}|>|z_{1}-z_{2}|>0$ are {\it analytic extensions of each other}.
We shall use also the similar notations and terminology below for other similar expressions and equalities.

Let $(W, Y_{W}, L_{W}(1))$ be an irreducible $V$-module. Assume that $W$ is graded either by $\Z$ or by $\Z+\frac{1}{2}$.
Also assume that for any irreducible $V$-module $(M, Y_{M}, L_{M}(1))$, the fusion rule $N_{WW}^{M}=1$ when $M$
is equivalent to $V$ and $N_{WW}^{M}=0$
when $M$ is not equivalent to $V$.
Assume in addition that there exist nondegenerate symmetric invariant bilinear forms   $(\cdot, \cdot)_{V}$ and $(\cdot, \cdot)_{W}$
on $V$ and $W$, respectively.

Let $V_{e}=V\oplus W$. We shall identify $V$ and $W$ with the corresponding subspaces
of $V_{e}$ so that $V$ and $W$ become subspaces of $V_{e}$. We define a vertex operator map
\begin{eqnarray*}
Y_{V_{e}}: \C^{\times}&\to &\hom(V_{e}\otimes V_{e}, \overline{V}_{e})\nn
z&\mapsto& Y_{V_{e}}(\cdot, z)\cdot
\end{eqnarray*}
by
\begin{eqnarray*}
Y_{V_{e}}(u, z)v&=&Y_{V}(u, z)v,\;\;\;u, v\in V,\\
Y_{V_{e}}(u, z)w&=&Y_{W}(u, z)w,\;\;\;u\in V, w\in W,\\
Y_{V_{e}}(w, z)v&=&Y_{WV}^{W}(w, z)v,\;\;\;v\in V, w\in W,\\
Y_{V_{e}}(w_{1}, z)w_{2}&=&Y_{WW}^{V}(w_{1}, z)w_{2},\;\;\;w_{1}, w_{2}\in W,
\end{eqnarray*}
where $Y_{WV}^{W}$ and $Y_{WW}^{V}$ are intertwining operators
of types ${W\choose WV}$ and ${V\choose WW}$, respectively, constructed
in \cite{FHL}. We recall their definitions:
$$Y_{WV}^{W}(w, z)v=e^{zL_{V}(-1)}Y_{W}(v, -z)w$$
for $v\in V$ and $w\in W$ and
\begin{equation}\label{ywwv}
(v, Y_{WW}^{V}(w_{1}, z)w_{2})_{V}=(Y_{WV}^{W}(e^{zL_{W}(1)}
e^{\pi i L_{W}(0)}z^{-2L_{W}(0)}w_{1}, z^{-1})v, w_{2})_{W},
\end{equation}
for $v\in V$ and $w_{1}, w_{2}\in W$, where $i=\sqrt{-1}$.
Note that in (\ref{ywwv}), since $W$ is graded by $\Z$ or $\Z+\frac{1}{2}$,
$z^{-2L_{W}(0)}w_{1}$ is well defined and involves only integer powers of $z$.
From (\ref{ywwv}), we also have
\begin{equation}\label{ywwv-1}
(w_{2}, Y_{WV}^{W}(w_{1}, z)v)_{W}
=(Y_{WW}^{V}(z^{-2L_{W}(0)}e^{-\pi i L_{W}(0)}e^{-z^{-1}L_{W}(1)} w_{1}, z^{-1})w_{2}, v)_{V}.
\end{equation}
Since $Y_{WV}^{W}$ and $Y_{WW}^{V}$ are intertwining operators,
they satisfy the $L(0)$-conjugation formula.

Let $\one_{V_{e}}=\one_{V}$. Let $L_{V_{e}}(0), L_{V_{e}}(-1), L_{V_{e}}(1): V_{e}\to V_{e}$
be the operators
that act as $L_{V}(0), L_{V}(-1)$, $L_{V}(1)$, respectively, on $V$ and as
$L_{W}(0), L_{W}(-1), L_{W}(1)$, respectively, on $W$.

Here is our second construction theorem:

\begin{thm}\label{second-const}
The triple $(V_{e}, Y_{V_{e}}, \one_{V_{e}}, L_{V_{e}}(1))$ is a
quasi-vertex operator algebra if $W$ is graded by $\Z$ and
is a quasi-vertex operator superalgebra if $W$ is graded by $\frac{1}{2}+\Z$. In the case that $V$ has a conformal element
$\omega_{V}$, $V_{e}$ also has a conformal element $\omega_{V_{e}}=\omega_{V}$.
\end{thm}
\pf
We first prove the skew-symmetry.  For $u, v\in V$ and $w\in
W$,  the skew-symmetry
$$Y_{V_{e}}(u, z)v=e^{zL_{V_{e}}(-1)}Y_{V_{e}}(v, -z)u$$
and
$$Y_{V_{e}}(u, z)w=e^{zL_{V_{e}}(-1)}Y_{V_{e}}(w, -z)u$$
follow from either the skew-symmetry for the
vertex operator algebra $V$ or the definition of $Y_{WV}^{W}(w, z)u$.

We now prove the skew-symmetry
\begin{equation}\label{skew-sym0}
Y_{V_{e}}(w_{1}, z)w_{2}=\epsilon_{W} e^{zL_{V_{e}}(-1)}Y_{V_{e}}(w_{2}, -z)w_{1}
\end{equation}
for $w_{1}, w_{2}\in W$, where $\epsilon_{W}$ is $1$ if $W$ is graded by $\Z$ and is $-1$ if $W$ is graded by $\Z+\frac{1}{2}$.
By definition, for $u\in V$ and $w_{1}, w_{2}\in W$,
\begin{eqnarray}\label{skew-sym}
\lefteqn{(u, Y_{V_{e}}(w_{1}, z)w_{2})_{V}}\nn
&&=(u, Y_{WW}^{V}(w_{1}, z)w_{2})_{V}\nn
&&=(Y_{WV}^{W}(e^{zL_{W}(1)}e^{\pi iL_{W}(0)}z^{-2L_{W}(0)}w_{1}, z^{-1})u, w_{2})_{W}\nn
&&=(e^{z^{-1}L_{W}(-1)}Y_{W}(u, -z^{-1})e^{zL_{W}(1)}e^{\pi iL_{W}(0)}z^{-2L_{W}(0)}w_{1}, w_{2})_{W}.
\end{eqnarray}
Using the fact that
the adjoint operator of $L_{W}(n)$ is $L_{W}(-n)$ for $n=-1, 0, 1$,
the  invariance of the bilinear form on $V$ and  the definitions
of $Y_{WV}^{W}$ and $Y_{WW}^{V}$, (\ref{ywwv-1}),
the right-hand side of (\ref{skew-sym}) is equal to
\begin{eqnarray}\label{skew-sym2}
\lefteqn{(e^{\pi iL_{W}(0)}z^{-2L_{W}(0)}w_{1}, e^{zL_{W}(-1)}Y_{W}
(e^{-z^{-1}L_{V}(1)}(-z^{2})^{L_{V}(0)}u, -z)e^{z^{-1}L_{W}(1)} w_{2})_{W}}\nn
&&=(e^{\pi iL_{W}(0)}z^{-2L_{W}(0)}w_{1}, Y_{WV}^{W}(e^{z^{-1}L_{W}(1)} w_{2}, z)
e^{-z^{-1}L_{V}(1)}(-z^{2})^{L_{V}(0)}u)_{W}\nn
&&=((-z^{2})^{L_{V}(0)}e^{-z^{-1}L_{V}(-1)}Y_{WW}^{V}(z^{-2L_{W}(0)}e^{-\pi i L_{W}(0)} w_{2},
z^{-1})e^{\pi iL_{W}(0)}z^{-2L_{W}(0)}w_{1}, u)_{V}.\nn
\end{eqnarray}
Since the weight of $L_{V}(-1)$ is $1$,  we have
\begin{eqnarray*}
(-z^{2})^{L_{V}(0)}e^{-z^{-1}L_{V}(-1)}&=&e^{zL_{V}(-1)}(-z^{2})^{L_{V}(0)}\nn
&=&e^{zL_{V}(-1)}z^{2L_{V}(0)}e^{-\pi iL_{V}(0)}.
\end{eqnarray*}
Using this formulas,  the $L(0)$-conjugation formula for $Y_{WW}^{V}$ and the fact that $e^{-2\pi iL_{W}(0)}=\epsilon_{W}$,
we see that the right-hand side of
(\ref{skew-sym2}) is equal to
\begin{eqnarray}\label{skew-sym3}
\lefteqn{(\epsilon_{W}e^{zL_{V}(-1)}Y_{WW}^{V}(w_{2}, -z)w_{1}, u)_{V}}\nn
&&=(\epsilon_{W}e^{zL_{V_{e}}(-1)}Y_{V_{e}}(e^{-2\pi iL_{W}(0)}w_{2}, -z)w_{1},u)_{V}\nn
&&=(u, \epsilon_{W}e^{zL_{V_{e}}(-1)}Y_{V_{e}}(e^{-2\pi iL_{W}(0)}w_{2}, -z)w_{1})_{V},
\end{eqnarray}
proving the skew-symmetry (\ref{skew-sym0}) and thus the skew-symmetry for $Y_{V_{e}}$ holds.

Next we prove the associativity. The associativity properties for the vertex operator maps $Y_{V}$ and $Y_{W}$
give the associativity
$$Y_{V_{e}}(u, z_{1})Y_{V_{e}}(v, z_{2})=Y_{V_{e}}(Y_{V_{e}}(u, z_{1}-z_{2})v, z_{2})$$
for $u, v\in V$ when $|z_{1}|>|z_{2}|>|z_{1}-z_{2}|>0$.
From \cite{FHL}, we obtain the following associativity properties
\begin{eqnarray*}
Y_{V_{e}}(u, z_{1})Y_{V_{e}}(w_{1}, z_{2})w_{2}
&=&Y_{V_{e}}(Y_{V_{e}}(u, z_{1}-z_{2})w_{1}, z_{2})w_{2},\\
Y_{V_{e}}(w_{1}, z_{1})Y_{V_{e}}(u, z_{2})w_{2}
&=&Y_{V_{e}}(Y_{V_{e}}(w_{1}, z_{1}-z_{2})u, z_{2})w_{2},\\
Y_{V_{e}}(w_{1}, z_{1})Y_{V_{e}}(w_{2}, z_{2})u
&=&Y_{V_{e}}(Y_{V_{e}}(w_{1}, z_{1}-z_{2})w_{2}, z_{2})u
\end{eqnarray*}
for $u\in V$ and $w_{1}, w_{2}\in W$ when $|z_{1}|>|z_{2}|>|z_{1}-z_{2}|>0$.

We still need to prove the associativity
$$Y_{V_{e}}(w_{1}, z_{1})Y_{V_{e}}(w_{2}, z_{2})w_{3}
=Y_{V_{e}}(Y_{V_{e}}(w_{1}, z_{1}-z_{2})w_{2}, z_{2})w_{3}$$
for $w_{1}, w_{2}, w_{3}\in W$ when $|z_{1}|>|z_{2}|>|z_{1}-z_{2}|>0$. By definition,
$$Y_{V_{e}}(w_{1}, z_{1})Y_{V_{e}}(w_{2}, z_{2})w_{3}
=Y_{WV}^{W}(w_{1}, z_{1})Y_{WW}^{V}(w_{2}, z_{2})w_{3}.$$
By assumption, the associativity property for intertwining operators for the quasi-vertex operator
algebra $V$  holds. In particular, there exist a $V$-module $M$ and
intertwining operators $\mathcal{Y}_{1}$ and $\mathcal{Y}_{2}$ of types ${W \choose MW}$
and ${M \choose WW}$, respectively,
such that
$$Y_{WV}^{W}(w_{1}, z_{1})Y_{WW}^{V}(w_{2}, z_{2})w_{3}
=\mathcal{Y}_{1}(\mathcal{Y}_{2}(w_{1}, z_{1}-z_{2})w_{2}, z_{2})w_{3}$$
when $|z_{1}|>|z_{2}|>|z_{1}-z_{2}|>0$. By assumption,  $M$ is equivalent to a direct sum of irreducible
$V$-modules  and the fusion rule $N_{WW}^{M}=1$ when $M$
is equivalent to $V$ and $N_{WW}^{M}=0$ when $M$ is not equivalent to $V$.
Then we see that $M$ must be equivalent to $V$ and $\mathcal{Y}_{2}$ is proportional to
$Y_{WW}^{V}$ after we identify $M$ with $V$.
On the other hand, it is a general fact that the fusion rule
$N_{VW}^{W}$ is $1$ (see
\cite{FHL}). In particular, $\mathcal{Y}_{1}$ is proportional to
$Y_{VW}^{W}(w_{1}, z_{1})$. Thus there exists $\lambda\in \C$ such that
\begin{equation}\label{assoc}
Y_{WV}^{W}(w_{1}, z_{1})Y_{WW}^{V}(w_{2}, z_{2})w_{3}
=\lambda Y_{VW}^{W}(Y_{WW}^{V}(w_{1}, z_{1}-z_{2})w_{2}, z_{2})w_{3}
\end{equation}
when $|z_{1}|>|z_{2}|>|z_{1}-z_{2}|>0$. From the genus-zero Moore-Seiberg equations (see \cite{MS} and \cite{H3}),
one can show that $\lambda=1$  and thus the associativity
property is proved.  Since the proof in this case is
easy, to make our proof self contained, here we give a direct proof that $\lambda=1$.

Note that
we have proved the skew-symmetry
\begin{equation}\label{skew-sym-3}
Y_{WW}^{V}(w_{1}, z)w_{2}=\epsilon_{W} e^{zL(-1)}Y_{WW}^{V}(w_{2}, -z)w_{1}
\end{equation}
for $w_{1}, w_{2}\in W$, where $\epsilon_{W}$ is $+1$ if $W$ is graded by $\Z$ and is $-1$ if $W$ is graded by $\Z+\frac{1}{2}$.
We shall use $\sim$ to mean that two expressions are analytic extensions of each other.
Then  for $w_{1}, w_{2}, w_{3}\in W$, using the associativity (\ref{assoc}) and the
definition of $Y_{WV}^{W}$, we have
\begin{eqnarray}\label{assoc-1}
\lefteqn{Y_{WV}^{W}(w_{1}, z_{1})Y_{WW}^{V}(w_{2}, z_{2})w_{3}}\nn
&&\sim\lambda Y_{W}(Y_{WW}^{V}(w_{1}, z_{1}-z_{2})w_{2}, z_{2})w_{3}\nn
&&\sim\lambda e^{z_{2}L(-1)}Y_{WV}^{W}(w_{3}, -z_{2})Y_{WW}^{V}(w_{1}, z_{1}-z_{2})w_{2}\nn
&&\sim\lambda^{2} e^{z_{2}L(-1)}Y_{W}(Y_{WW}^{V}(w_{3}, -z_{1})w_{1}, z_{1}-z_{2})w_{2}.
\end{eqnarray}
On the other hand, using (\ref{assoc}) and (\ref{skew-sym-3}), we have
\begin{eqnarray}\label{assoc-2}
\lefteqn{Y_{WV}^{W}(w_{1}, z_{1})Y_{WW}^{V}(w_{2}, z_{2})w_{3}}\nn
&&\sim \epsilon_{W} Y_{WV}^{W}(w_{1}, z_{1})e^{z_{2}L(-1)}Y_{WW}^{V}(w_{3}, -z_{2})w_{2}\nn
&&\sim \epsilon_{W}e^{z_{2}L(-1)}Y_{WV}^{W}(w_{1}, z_{1}-z_{2})Y_{WW}^{V}(w_{3}, -z_{2})w_{2}\nn
&&\sim \epsilon_{W}\lambda e^{z_{2}L(-1)}Y_{W}(Y_{WW}^{V}(w_{1}, z_{1})w_{3}, -z_{2})w_{2}\nn
&&\sim \epsilon_{W}^{2}\lambda e^{z_{2}L(-1)}Y_{W}(e^{z_{1}L(-1)}Y_{WW}^{V}(w_{3}, -z_{1})w_{1}, -z_{2})w_{2}\nn
&&\sim \lambda e^{z_{2}L(-1)}Y_{W}(Y_{WW}^{V}(w_{3}, -z_{1})w_{1}, z_{1}-z_{2})w_{2}.
\end{eqnarray}
Since both the right-hand sides of (\ref{assoc-1}) and (\ref{assoc-2}) are
defined on the same region and they are analytic extensions of each other, they
must be equal on the region that they are defined. Thus $\lambda=1$.

Note that the skew-symmetry (\ref{skew-sym0}) has a sign $\epsilon_{W}$.
Since $Y_{V_{e}}$ satisfies the associativity and skew-symmetry,
from \cite{H3}, it also satisfies the commutativity for quasi-vertex operator algebras
when $W$ is graded by $\Z$ and satisfies the commutativity for quasi-vertex operator
superalgebras when $W$ is graded by $\Z+\frac{1}{2}$. Since $Y_{V_{e}}$ involves
only integer powers of the variable,  the rationality for products and iterates also hold (see also \cite{H3}).
The other axioms for grading-restricted vertex operator (super)algebras can be easily verified.
Thus we have proved that
$(V_{V_{e}}, Y_{V_{e}}, \one_{V_{e}}, L_{V_{e}}(1))$
is a quasi-vertex operator algebra when $W$ is graded by $\Z$ and
is a quasi-vertex operator superalgebra when $W$ is graded by $\Z+\frac{1}{2}$.
\epfv

\begin{exam}
{\rm As is mentioned above, this second construction is in fact a generalization of
the construction of the moonshine
module vertex operator algebra and the corresponding vertex operator superalgebra in \cite{H1}.
In fact, take $V$ to be the
fixed point vertex operator subalgebra $V_{\Lambda}^{+}$
of the Leech lattice vertex operator algebra $V_{\Lambda}$ under the automorphism induced from
the automorphism $\alpha\mapsto -\alpha$ of the Leech lattice $\Lambda$. Take $W$ to be the
fixed point $V_{\Lambda}^{+}$-submodule $(V_{\Lambda}^{T})^{+}$ of the irreducible twisted
$V_{\Lambda}$-module $V^{T}_{\Lambda}$. Then the conditions to use Theorem \ref{second-const}  are
satisfied. By Theorem \ref{second-const} and the fact that  $V_{\Lambda}^{+}$ has a conformal element, the moonshine module
$V^{\natural}=V_{\Lambda}^{+}\oplus (V_{\Lambda}^{T})^{+}$
has a structure of vertex operator algebra for which the vertex operator map, the vacuum
and the conformal element are given above.
If we take $V$ to be the same and $W$ to be the eigenspace $(V_{\Lambda}^{T})^{-}$
for the action of the automorphism above on $V^{T}_{\Lambda}$
with the eigenvalue $-1$. Then $V_{\Lambda}^{+}\oplus (V_{\Lambda}^{T})^{-}$
has a structure of vertex operator superalgebra given above. See \cite{FLM},
\cite{H1} and \cite{H2} for the background  material that can be used to verify the conditions to use
Theorem \ref{second-const}. }
\end{exam}

\begin{exam}
{\rm Another application is a construction of the vertex  operator
superalgebra structure on the direct sum $L(\frac{1}{2}, 0)\oplus
L(\frac{1}{2}, \frac{1}{2})$ of the vertex operator algebra
$L(\frac{1}{2}, 0)$ for the minimal
model of central charge $\frac{1}{2}$ and the $L(\frac{1}{2}, 0)$-module
$L(\frac{1}{2}, \frac{1}{2})$ of lowest weight $\frac{1}{2}$. The conditions to use
Theorem \ref{second-const} are satisfied by the results obtained in \cite{W} and
\cite{H3.1} or \cite{H4}. Thus $L(\frac{1}{2}, 0)\oplus
L(\frac{1}{2}, \frac{1}{2})$ has a structure of vertex operator superalgebra given above.
After tensoring the vertex operator algebra
$L(\frac{1}{2}, 0)\oplus
L(\frac{1}{2}, \frac{1}{2})$ with the vertex operator algebra associated to the lattice
of rank $1$ generated by $\alpha$ satisfying $(\alpha, \alpha)=m\in 2\Z_{+}$,
we obtain the vertex operator superalgebra for the Moore-Read state with
the filling factors $\nu=\frac{1}{m}$
in the conformal-field-theoretic study of fractional quantum Hall states.
We refer the reader to \cite{H11} for details on these examples.}
\end{exam}

\noindent {\small \sc Department of Mathematics, Rutgers University,
110 Frelinghuysen Rd., Piscataway, NJ 08854-8019}
\vspace{1em}

\noindent  {\it and}
\vspace{1em}

\noindent {\small \sc Beijing International Center for Mathematical
Research, Peking University, Beijing 100871, China}
\vspace{1em}

\noindent {\em E-mail address}: yzhuang@math.rutgers.edu

\end{document}